\documentclass[11pt]{article}

\usepackage[english]{babel}
\usepackage{mathtools}
\usepackage{amsmath,amssymb,dsfont}
\usepackage{epsfig, subfigure}
\usepackage{natbib}
\usepackage{color}

\usepackage{wrapfig} 
\usepackage{bbm}

\usepackage{appendix}
\usepackage{versions}
\usepackage{includex} 

\usepackage{fancyhdr}
\usepackage[colorlinks,citecolor=blue,urlcolor=blue]{hyperref}
\usepackage{hypernat}

\usepackage{url}

\usepackage{multirow}
\usepackage{booktabs}
\usepackage{slashbox}

\let\proglang=\textsf
\newcommand{\pkg}[1]{{\fontseries{b}\selectfont #1}}

\DeclareMathOperator{\supp}{supp}

\newcommand{\lp}{\left(} 
  \newcommand{\rp}{\right)}
\newcommand{\lb}{\left\{} 
  \newcommand{\rb}{\right\}}

\usepackage{mathabx} 
\newcommand{\piest}{\widecheck{\pi}_{n}}
\newcommand{\uestI}{\widecheck{u}_{1,n}}
\newcommand{\uestII}{\widecheck{u}_{2,n}}
\newcommand{\gestnh}{\widecheck{g}_{n,h}}

\newcommand{\PP}{{\mathbb P}}

\newcommand{\GG}{{\mathbb G}}
\newcommand{\GGn}{{\mathbb G}_n}

\newcommand{\RR}{\mathbb{R}}

\newcommand{\argmax}{\operatorname{argmax}}

\newtheorem{theorem}{Theorem}[section]
\newtheorem{prop}{Proposition}[section]
\newtheorem{corollary}{Corollary}[section]
\newtheorem{lemma}{Lemma}[section]

\usepackage[usenames,dvipsnames]{xcolor}

\excludeversion{supplement-file}
\excludeversion{main-file}
\includeversion{all-in-one-file}

\begin{supplement-file}
\title{Supplementary Material for ``Inference for a two-component mixture
    of symmetric distributions under log-concavity''}
\end{supplement-file}

\begin{main-file}
  \title{Inference for a two-component mixture
    of symmetric distributions under log-concavity}
\end{main-file}

\begin{all-in-one-file}
  \title{Inference for a two-component mixture of symmetric distributions
    under log-concavity}
\end{all-in-one-file}

\author{Fadoua Balabdaoui and Charles R.\ Doss}
\date{}

\begin{document}
\maketitle

\begin{all-in-one-file}
  \begin{abstract}

    In this article, we revisit the problem of estimating the unknown
    zero-symmetric distribution in a two-component location mixture model,
    considered in previous works, now under the assumption that the
    zero-symmetric distribution has a log-concave density.  When consistent
    estimators for the shift locations and mixing probability are used, we
    show that the nonparametric log-concave Maximum Likelihood estimator
    (MLE) of both the mixed density and that of the unknown zero-symmetric
    component are consistent in the Hellinger distance. In case the
    estimators for the shift locations and mixing probability are $\sqrt
    n$-consistent, we establish that these MLE's converge to the truth at the
    rate $n^{-2/5}$ in the $L_1$ distance.  To estimate the shift locations
    and mixing probability, we use the estimators proposed by
    \cite{hunteretal2007}. The unknown zero-symmetric density is efficiently
    computed using the \proglang{R} package \pkg{logcondens.mode}.

  \end{abstract}
\end{all-in-one-file}

\begin{main-file}
\begin{abstract}

  In this article, we revisit the problem of estimating the unknown zero-symmetric
  distribution in a two-component location mixture model, considered in
  previous works, now under the assumption that the zero-symmetric distribution
  has a log-concave density.  When consistent estimators for the shift
  locations and mixing probability are used, we show that the nonparametric
  log-concave Maximum Likelihood estimator (MLE) of both the mixed density
  and that of the unknown zero-symmetric component are consistent in the Hellinger
  distance. In case the estimators for the shift locations and mixing
  probability are $\sqrt n$-consistent, we establish that these MLE's
  converge to the truth at the rate $n^{-2/5}$ in the $L_1$ distance.  To
  estimate the shift locations and mixing probability, we use the estimators
  proposed by \cite{hunteretal2007}. The unknown zero-symmetric density is
  efficiently computed using the \proglang{R} package \pkg{logcondens.mode}.

\end{abstract}
  \end{main-file}

\begin{supplement-file}
\begin{abstract}
  This supplement contains proofs and other technical material for
  ``Inference for a mixture of symmetric distributions under log-concavity''.
  Reference numbers to equations, theorems, or other statements in that
  document have a unified numbering system with this document.  The
  bibliography is shared and found at the end of the main document.
\end{abstract}
\end{supplement-file}

\section{Introduction}

Let us assume that $X_1,\ldots, X_n$ are independent and identically
distributed (i.i.d.) draws from a mixture distribution, with cumulative
distribution function (cdf) $G^0$ given by
\begin{equation}
  \label{eq:general-mixture}
  G^0(x) = \sum_{i=1}^k \pi_i^0 F_i^0(x), \;\;\;\; x \in \RR,
\end{equation}
for some integer $k \ge 2$, where $F_i^0$ are cdfs, $\pi_i^0 \ge 0,$ and
$\sum_{i=1}^k \pi_i^0 = 1$.  Such mixture distributions are very common in
statistical modeling, in part because a variety of data generating frameworks
lead to mixture models; for instance, one common approach to clustering
problems leads to estimation of a mixture density \citep{Fraley:2002bg}.
Another reason for this popularity is that they are very flexible and many
distributions can be well approximated by some mixture model (see, e.g.,
\cite{Everitt:1981ua}, \cite{Titterington:1985vh}, or
\cite{McLachlan:2000ga}).

In this paper, we revisit the semi-parametric mixture model already studied
by \cite{bordes2006} and \cite{hunteretal2007}. In this model, it is assumed
that the mixing distributions $F_i, 1 \le i \le k$ in
\eqref{eq:general-mixture} are such that
\begin{eqnarray*}
  F_i(x) = F^0(x - u_i^0)
\end{eqnarray*}
for $u_i^0 \in \RR$, $i=1,\ldots, k$, and $F^0$ is a distribution function restricted to be
symmetric about $0$, i.e.\  $F^0(-x) = 1-F^0(x-)$.
This model was also studied more  recently by \cite{ButuceaVand2014}.
All of these authors have actually focused on the case $k = 2$:
\begin{equation}
  \label{eq:MixtureModel}
  G^0(x) =  \pi^0 F^0(x-u_1^0) + (1- \pi^0) F^0(x- u_2^0), \;\;\;\; x \in \RR.
\end{equation}
This is still, in fact, a flexible model which is useful in many scenarios
(see our data applications in Section~\ref{sec:data}).  As the main goal is
to estimate the mixing parameters and the mixing component $F^0$, it is
crucial to be assured that there exists a unique solution $(\pi^0, u^0_1,
u^0_2, F^0)$ for a given $G^0$ determined by \eqref{eq:MixtureModel}.
\cite{bordes2006} and \cite{hunteretal2007} were able to establish that
identifiability holds under some suitable conditions on the mixing
parameters. Their result states that if $u_1^0 < u_2^0$ and $\pi^0 \notin
\left\{ 0, 1/2, 1 \right\}$, then $G^0$ given \eqref{eq:MixtureModel} is
identifiable for any zero-symmetric distributions $F^0$. Furthermore, the
condition is necessary and sufficient since any distribution $G^0$ that is
symmetric about its median clearly cannot be $2$-identifiable; see Theorem 2
of \cite{hunteretal2007}.

After having shown identifiability, \cite{hunteretal2007} put their focus on
estimating $(\pi^0,u_1^0,u_2^0, F^0)$.  They have shown that their estimator
of the parametric component $(\pi^0,u_1^0,u_2^0)$ is consistent and
asymptotically normal.  However, the obtained estimator of $F^0$ is not even
guaranteed to have the properties of a genuine cdf (i.e., it is not
necessarily nondecreasing).  On the other hand, \cite{bordes2006},
\cite{ButuceaVand2014}, and \cite{cheeWang13} use a KDE approach to
estimation of $F^0$. The resulting estimators are proper distribution
functions, but the procedures involve a model-selection procedure
(cross-validation or Akaike or Bayesian information criterion) to choose the
tuning parameter.  The estimators of \cite{hunteretal2007} and
\cite{ButuceaVand2014} for the mixture parameters are shown to converge
weakly to a multivariate Gaussian at the parametric rate $n^{-1/2}$ under
some regularity conditions on $F^0$ which are related to smoothness in the
case of \cite{ButuceaVand2014}.  \cite{bordes2006} obtain also a convergence
rate under smoothness assumptions, but their rate of convergence is much
slower (of order $n^{-1/4 + \alpha}$, for any $\alpha >
0$). \cite{bordes2006} show that the same rates of convergence are inherited
by their kernel estimator of $F^0$ in the supremum norm, under the assumption
that the location parameters $u^0_1$ and $u^0_2$ are unknown.
If $F^0$ is assumed to admit a density $f^0$, then
\cite{bordes2006} provide only almost sure
consistency in the supremum norm.
For their kernel estimator, \cite{ButuceaVand2014}) obtain, for pointwise convergence, a rate of order
$n^{-(2\beta-1)/(4\beta)}$ in the quadratic risk assuming smoothness of level
$\beta > 1/2$ and assuming that the bandwidth is chosen optimally
(the authors suggest using cross validation).

Hence, the proposed estimators of $F^0$ in the aforementioned works suffer
various practical difficulties, including slow rates of convergence (or
as-of-yet unknown rates) for it or its density, the estimator not being a proper cdf, or the need for model-selection procedures to choose a tuning
parameter. Our goal in this paper is to circumvent those issues by
constructing an estimator of the density $f^0$ which
\begin{itemize}
\item converges to the truth with a provably good convergence rate,
\item can be efficiently computed,
\item does not require a tuning parameter,  and furthermore,
\item is unimodal.
\end{itemize}
Unimodality is a natural constraint to enforce; when using a mixture model,
it is somewhat unnatural to imagine a multimodal mixture component
density.
However, using unimodality involves some technical difficulties: enforcing
unimodality on $f^0$ is not directly feasible, because the class of unimodal
densities is too large and the MLE of a unimodal density fails to exist even
in the simple one-dimensional setting (with no mixing). We propose instead to
assume that $f^0$ satisfies the shape constraint of log-concavity (i.e.,
$\log f^0$ is
concave).

Log-concave functions are always unimodal, and have been used to great
success in nonparametric modeling.  Unlike the class of unimodal densities,
the log-concave class admits an MLE \citep{walther_02}.  Many papers have
studied the log-concave MLE on $\RR$ or $\RR^d$ and much is already known about its
large sample properties, both local and global; see e.g.  \cite{pal_07},
\cite{duembgen_09}, \cite{balabdaoui_09}, \cite{cule_10_ejs}, \cite{cule_10},
\cite{duembgen_10}, \cite{schuhmacher_10}, \cite{chen_2013},
\cite{dosswellner13}, and \cite{Kim:2014wa}.  \cite{balabdaoui_13} studied
asymptotics and confidence intervals of the discrete log-concave MLE of a
probability mass function in the well- and misspecified settings.
\cite{duembgen_07}, \cite{rufibach_07} and \cite{duembgen_logcon10} study
algorithms for computation of the MLE, allowing unequal weights to be
assigned to the observations, an important feature of which we will take
advantage.

In the present context, we need to consider the class of {\em zero-symmetric}
log-concave densities on $\RR$, which has not been considered before.  To do
so, we note that if $f$ is zero-symmetric and log-concave on $\RR$, then $f^+(t)
:= 2 f(t) \mathbb{I}_{t \in [0, \infty)}$ is log-concave with mode at 0.
Thus, through a simple transformation of the data, it can be shown that the
original estimation problem is equivalent to maximizing the log-likelihood
over the class of log-concave densities on $[0,\infty)$ with mode at 0.  We
can then compute the maximum of the log-likelihood easily by alternating
between the EM algorithm \citep{MR0501537} and the active set algorithm
provided in the \proglang{R} package \pkg{logcondens.mode} which computes the
log-concave MLE with a fixed mode.  We use the fact that the active set
algorithm allows for unequal weights to be assigned to the data points: here,
the weights assigned are proportional to the posterior probabilities from the
EM algorithm.

We are able to show that the zero-symmetric log-concave MLE converges
in probability
to the true zero-symmetric log-concave component density in the Hellinger
distance and in the supremum norm on sets of continuity of the true
density. Furthermore, it can be
shown that our estimator converges to the truth at the rate $n^{-2/5}$ in the
$L_1$-distance.
Although the risk measure we use here is different from the one considered by
\cite{ButuceaVand2014}, it seems that the rate of convergence of our MLE,
when the true mixture component is log-concave, is faster than that given in
their Theorem~4 for their KDE when the smoothness parameter $\beta$ satisfies
$\beta \in (1/2, 5/2)$.  Note for an estimator $\hat g_n$ of $g_0$ in the
direct density estimation problem based on i.i.d.\ observations from $g_0$
(as opposed to the mixture setting) when $g_0$ has smoothness $\beta$ the
optimal pointwise rate of convergence of $| \hat g_n(x_0) - g_0(x_0)|$ at a
fixed point $x_0$ is $ n^{- \beta / (2 \beta + 1)}$
\citep{Stone:1980vq}.  Note also that \cite{duembgen_09} find a rate of
convergence of $ (\log n)^{\frac{\beta}{2 \beta + 1}} n^{-\frac{\beta}{2
    \beta + 1}}$ in the uniform norm on compact sets for the log-concave MLE,
in the direct density estimation problem, when the true density $g_0$ is
log-concave and also lies in a H\"older class with smoothness $\beta$, i.e.\
\begin{eqnarray*}
  \vert g_0(x) - g_0(y) \vert &\le& L \vert x - y \vert, \ \ \textrm{if $\beta =1$} \\
  \vert g_0'(x) - g_0'(y) \vert &\le & L \vert x-y \vert^{\beta-1}, \ \ \textrm{if $\beta > 1$}
\end{eqnarray*}
for some $L > 0$.
This rate is optimal for nonparametric estimation with smoothness $\beta$
(the log
factor being due to the supremum norm \citep{khasminskii1978}), and no
bandwidth needs to be chosen.

We note that, although we refer to our estimator as the log-concave MLE, we
do not use a \lq\lq pure\rq\rq \ maximum likelihood approach since we feed in
other estimators of $(\pi^0, u_1^0, u_2^0)$ to our likelihood, which we
maximize to estimate $f^0$ and thus $g^0$, the density of the mixed
distribution $G^0$.  An alternative approach is to estimate both the
parametric and nonparametric components simultaneously by maximum likelihood.
However, there are many additional difficulties in that approach, due
to the complicated non-concave nature of the log-likelihood function; see
Section~\ref{sec:finite-sample-estimator}.

We also note that we are not the first to use log-concavity in mixture
modeling; \cite{Chang:2007ki} and \cite{Eilers:2007hc} consider univariate
mixtures of log-concave densities, and \cite{cule_10} consider multivariate
mixtures of log-concave densities.  However, in none of those settings was
symmetry imposed, perhaps because the authors were not worried about the
(often fundamental) question of identifiability.  Thus, their work does not
directly apply in our setting.

The paper will be structured as follows.
In %
Section~\ref{sec:finite-sample-estimator} we establish existence of the MLE
and provide a necessary condition for a candidate to be equal to the
estimator.
In Section~\ref{consistency},
we establish
consistency in the Hellinger distance. This implies other forms of
consistency by the results of \cite{cule_10_ejs}.  The techniques we used are
re-adapted from \cite{pal_07}, \cite{cule_10_ejs} and
\cite{schuhmacher_10} to deal with the
additional difficulties of a mixture model.
In Section~\ref{sec:rates},
we find that
the MLEs of $f^0$ and $g^0$ converge to the truth at a rate of order
$n^{-2/5}$.
In Section~\ref{NumSection}, we develop a likelihood ratio procedure based on our estimator in the problem
of testing absence of mixing.  We also consider the problem of clustering where we use the estimators of the posterior probabilities obtained via our log-concave MLE.  In both problems, we compare our method to alternative or existing approaches.
In Section~\ref{sec:data},
we present two data
applications. Section 7 gathers some conclusions.
Proofs and technical details can be found in the online supplementary material.

\section{The model and estimation via Maximum likelihood}
\label{sec:finite-sample-estimator}

Let $X_1, \ldots, X_n$ to be $n$ independent observations assumed to come from
the location mixture with cdf $G^0$ which we now assume has a density, given by
\begin{eqnarray}
  \label{eq:defnmix}
  g^0(x) = \pi^0 f^0(x- u^0_1) + (1- \pi^0) f^0(x-u^0_2),
\end{eqnarray}
for some $\pi^0 \in (0,1) \setminus \{1/2\}$, $u^0_1, u^0_2 \in \mathbb R$
such that $u^0_1 \ne u^0_2$.
We assume that $f^0$ is a zero-symmetric log-concave density, i.e.\ $f^0 \in
{\cal SLC}_1$ where
\begin{align*}
  {\cal SLC}_1 := {\cal SLC} \cap \lb f: \int_{\RR} f(u)du = 1 \rb,
  \; \mbox{and} \;
  {\cal SLC}  := \lb e^{\psi} : \psi \in {\cal SC} \rb  ,
\end{align*}
and ${\cal SC}$ is the class of concave functions on $\RR$ that are upper
semi-continuous (``closed'') and proper \citep{rock70}, and satisfy
$\psi(x)=\psi(-x)$.
The upper semi-continuity condition is made only for the purpose
of uniqueness.  Then
\begin{eqnarray}\label{LL}
  L(\pi,u_1,u_2, f)
  := \sum_{j=1}^n \log\Big(\pi f(X_j - u_{1}) + (1- \pi) f(X_j - u_{2}) \Big).
\end{eqnarray}
is the log-likelihood in this problem.  In the case of estimation of a
log-concave density on $\RR$, the log-likelihood is a concave function
(\cite{pal_07}, \cite{rufibach_07}, \cite{duembgen_09}).  However,
\cite{Duembgen:2010-arxiv-v3} study a semiparametric model incorporating
log-concavity and find a non-concave likelihood; see their Section~3.3
including a plot on page 18.  Unfortunately, our objective function $L$ is
also far from concave.  Consider order statistics $X_{(1)}, \ldots, X_{(n)}$,
a fixed $\pi \in (0,1) \setminus \{1/2\}$, and (zero-symmetric log-concave)
$f$ with support given by $\supp(f) = [-s,s]$ and $\inf_{x \in \supp(f)} f(x)
> 0$.  Assume $u_1 < u_2$ are such that $[u_1-s, u_1+s] \cup [u_2-s, u_2+s]
\supset [X_{(1)}, X_{(n)}]$ so that $L(\pi,u_1,u_2,f) > -\infty.$ Let $j :=
\min \lb i : X_{(i)} > u_2-s \rb$ be the index of the smallest order
statistic contained in the support of the second component,
and let $\tilde \delta := X_{(j)} - (u_2-s)$.  Then not only
does $L(\pi, u_1, \cdot, f)$ fail to be concave, but it is in fact discontinuous
at $u_2 + \tilde \delta$.

\bigskip

\par \noindent

We now describe our estimation approach.  Let  $\piest, \uestI,\uestII$ be estimators of
$\pi, u_1,u_2$, where we assume $\piest \in (0,1) \setminus \lb 1/2 \rb$
and $\uestI < \uestII$. We will generally think of these estimators as being
$\sqrt{n}$-consistent.
We will then consider
maximizing the
log-likelihood
\begin{eqnarray}\label{LL}
  f \mapsto \sum_{j=1}^n \log\Big(\piest f(X_j - \uestI) + (1- \piest) f(X_j - \uestII) \Big)
\end{eqnarray}
over $\mathcal{SLC}_1$.
Using the Lagrange penalty term introduced by \cite{silverman82}, this is
equivalent to maximizing the criterion $\Phi_n$ defined as
\begin{eqnarray*}
  \Phi_n(\psi) = \frac{1}{n} \sum_{j=1}^n \log\Big[\piest e^{\psi(X_i-
    \uestI)} + (1- \piest) e^{\psi(X_i - \uestII)}\Big] - \int_{\mathbb R} e^{\psi(x)} dx
\end{eqnarray*}
over  $\mathcal{SC}$.
We will abusively use the term MLE for
our estimators of $f^0$ and $g^0$ despite the fact that
the mixing parameters $(\pi, u_1,u_2)$
are not a part of the space over which the likelihood is maximized.
In the next proposition we establish existence of the MLE, and describe its
nature.

\bigskip

\begin{prop}\label{Existence}
  The criterion $\Phi_n$ admits a maximizer $\widehat \psi_n$. Furthermore,
  the following holds true almost surely, letting $\widehat f_n =
  e^{\widehat \psi_n}$.
  \begin{itemize}
  \item $\widehat f_n$ is in $\mathcal{SLC}_1$.
  \item For $i =1, \ldots, n$, let
    \begin{equation}
      \label{eq:1}
      Z_{2i-1} = \vert X_i-\uestI \vert  \mbox{ and } Z_{2i} = \vert X_i - \uestII\vert.
    \end{equation}
    Then, on $[0, \infty)$ the MLE $\widehat \psi_n$ changes slope only at points belonging to the set
    \begin{eqnarray*}
      \big \{Z_1 ,  Z_2 , \ldots,  Z_{2n-1} ,  Z_{2n} \big \}.
    \end{eqnarray*}
    Furthermore, $\widehat \psi_n'(0) =0$, and  $\widehat \psi_n(x) = -\infty$
    if and only if  $x \notin [-Z_{(2n)}, Z_{(2n)}]$ where $Z_{(2n)}$ is the largest order
    statistic of $Z_1 , Z_2 , \ldots, Z_{2n-1} , Z_{2n} $.
  \end{itemize}

\end{prop}

\bigskip

The MLE of $f$ will be denoted by $\widehat f_n$ throughout, and that of
$g$ by $\widehat g_n$.  In the following, we give a
necessary condition for a log-concave function $f = \exp(\psi)$ to be the
MLE.

Proposition~\ref{prop:necessary-condition} is
interesting to compare with the characterization of \cite{duembgen_09} for the log-concave MLE.
The result is also useful in combination
with the EM-algorithm described below as its non-fulfillment indicates that
convergence is not yet reached.

\medskip

\begin{prop}
  \label{prop:necessary-condition}
  Let $\psi$ be a zero-symmetric concave function on $\mathbb R$ such that
  $\psi(x) = -\infty$ if and only if $x \notin [-Z_{(2n)}, Z _{(2n)}]$ where
  $Z_{(2n)}$ is defined in Proposition~\ref{Existence}, and $\psi'(0)=0$. If $\exp(\psi) =
  \widehat f_n$ is the MLE, then for any real zero-symmetric function $\Delta$
  such that $\psi + \epsilon \Delta \in \mathcal {SC}$ for some $\epsilon >
  0$ we have that
  \begin{eqnarray}\label{CharMLE}
    \frac{1}{n} \sum_{i=1}^n \Big \{\widehat p_n(X_i) \Delta(X_i - \uestI) +
    (1-\widehat p_n(X_i)) \Delta (X_i- \uestII) \Big \} \le \int_{\mathbb R} \widehat{f}_n(x) \Delta(x) dx
  \end{eqnarray}
  where
  \begin{eqnarray}\label{PostWeights}
    \widehat p_n(X_i) = \frac{\piest \widehat f_n(X_i- \uestI)}{\piest \widehat
      f_n(X_i- \uestI)
      + (1-\piest) \widehat{f}_n(X_i-\uestII)} = \frac{\piest \widehat
      f_n(X_i- \uestI)}{\widehat g_n(X_i)},
  \end{eqnarray}
  for $i =1, \ldots, n$.
\end{prop}
\medskip

Next, we give the condition in (\ref{CharMLE}) under an alternative form.
\cite{duembgen_09} shows that the log-concave MLE is uniquely characterized
by the fact that the first integral of the cdf of the MLE stays below the
first integral of the empirical distribution, while touching it exactly at
the points where the logarithm of the MLE changes slope. To derive a related
result, let $\widehat{\mathbb F}_n$ denote the cdf of the discrete
distribution putting mass $\widehat p_n(X_i)/n$ at $ Z_{2i-1} $ and $
(1-\widehat p_n(X_i))/n$ at $Z_{2i}$ for $i=1, \ldots, n$, where
$\widehat{p}_n(X_i)$ was defined in (\ref{PostWeights}) and $Z_i$ was defined
in \eqref{eq:1}.  That is,
\begin{eqnarray*}
  \widehat{\mathbb F}_n = \frac{1}{n} \sum_{i=1}^n \left(\widehat p_n(X_i) \delta_{Z_{2i-1}}  + (1-\widehat p_n(X_i)) \delta_{Z_{2i}}  \right)
\end{eqnarray*}
where $\delta_x(t) = 1_{\{[x, \infty)\}}(t)$.
Let $\widehat f^+_n(x)= 2\widehat f_n(x) 1_{x \in [0, \infty)}$, $\widehat
\psi_n^+ = \log(\widehat f_n^+)$ and let $\widehat F^+_n$ be the cdf of $\widehat
f^+_n$.

\medskip

\begin{prop}\label{AltChar}
  If $\widehat f_n$ is the MLE of the
  component $f^0 \in \mathcal{SLC}_1$ then
  \begin{eqnarray}\label{CharMLE2}
    \int_0^z \widehat{F}^+_n (x) dx
    \begin{cases}
      \le \int_0^z \widehat{\mathbb F}_n(x)dx, \ \textrm{for $z \in [0, Z_{(2n)}]$} \\
      = \int_0^z \widehat{\mathbb F}_n(x) dx, \ \textrm{if $\widehat \psi^+_n(z-) > \widehat \psi^+_n(z+)$}.
    \end{cases}
  \end{eqnarray}

\end{prop}

\medskip

\section{Consistency}
\label{consistency}

The main result of this section is to establish consistency in the Hellinger
distance of the MLEs $\widehat g_n$ and $\widehat f_n$ as $n \to \infty$,
where the Hellinger distance is defined by
$$H(p,q) := \sqrt{(1/2) \int \lp
  \sqrt{p(x)}-\sqrt{q(x)} \rp^2 dx }.$$  We will also find consistency for
$\widehat f_n$ in certain exponentially weighted metrics.  Our approach to
the problem follows the idea of \cite{pal_07} and \cite{duembgen_09} but will
require handling carefully the extra complexity induced by the mixture. As in
\cite{pal_07}, \cite{cule_10_ejs} and \cite{schuhmacher_10}, we will first
need to show that the MLE of the mixed density and hence the MLE of the
log-concave component are bounded. Here, the claimed boundedness will be only
in probability, which is weaker than the almost sure boundedness proved in
the aforementioned articles.  Those articles, however, were able to take
advantage of the fact that the level sets of a bounded unimodal function are
convex and compact; such a statement does not hold if we consider a mixture
of two unimodal functions instead of a single unimodal function, even if the
two components are log-concave.   So, instead of studying how the empirical
distribution
behaves over the class of compact intervals, we will instead need to
study its behavior over more complicated classes of functions.  This is what
is done in Propositions %
\ref{Glivenko1} and \ref{Glivenko2}. %

\bigskip

\begin{theorem}\label{ConvMixed}
  Let $g^0$ be as in \eqref{eq:defnmix} and $\widehat g_n$ be the MLE
  of $g^0$. Then we have that
  \begin{eqnarray*}
    H(\widehat g_n, g^0)  = o_p(1).
  \end{eqnarray*}
\end{theorem}

\medskip

\bigskip

\par \noindent Consistency of the log-concave component, $\widehat f_n$, follows now from Theorem \ref{ConvMixed}.
\medskip

\begin{corollary}
  \label{cor:consistency-f}
  Let $f^0$ denote again the true log-concave zero-symmetric density. Then,
  \begin{eqnarray*}
    H(\widehat f_n, f^0)= o_p(1),
  \end{eqnarray*}
  and for any $a \in (0, a_0)$ such that $f^0(x) \le \exp(-a_0 x + b) $ for some $b \in \RR$, then
  \begin{eqnarray*}
    \int_{\RR} e^{a t}  \vert \widehat f_n(t) - f^0(t) \vert = o_p(1),
  \end{eqnarray*}
  and
  \begin{eqnarray*}
    \sup_{t \in [-A,A]} e^{a t}  \vert \widehat f_n(t) - f^0(t) \vert = o_p(1).
  \end{eqnarray*}
  on any continuity set $[-A,A]$ of $f^0$, where $A$ may be $\infty$ if $f^0$
  is continuous on all of $\RR$.
\end{corollary}

\medskip

\section{Rates of convergence}
\label{sec:rates}
In this section, we aim at refining the convergence result obtained in the
previous section to attain a rate of convergence for both $\widehat f_n$ and
$\widehat g_n$ in the $L_1$ distance.  To this goal, we need first to recall
some definitions from empirical processes theory. Given a class of functions
$\mathcal F$, the bracketing number of $\mathcal F$ under some distance
$\Vert \cdot \Vert $ is defined as
\begin{eqnarray*}
  N_{[ \ ]}(\epsilon, \mathcal F, \Vert \cdot \Vert )= \min \Big \{k: \exists \ \underline{f}_1, \bar{f}_1, \ldots, \underline{f}_k, \bar{f}_k \ \textrm{s.t.}  \Vert \underline{f}_i - \bar{f}_j \Vert \le \epsilon, \mathcal F \subset  \cup_{i=1}^k [\underline{f}_i, \bar{f}_i] \Big\}
\end{eqnarray*}
where $[l, u] = \{f : \in \mathcal F: l \le f \le u \}$.   In this section, we refine the consistency result above by deriving the rate of convergence of the MLE's $\widehat g_n$ and $\widehat f_n$ of the mixed density and the zero-symmetric log-concave component respectively.

\medskip

For fixed $M > 0$, $a_0 < b_0$ and $\delta \in  \left(0, (b_0 - a_0)/2\right)$, consider the class of functions
\begin{eqnarray*}
  \mathcal{G}  & = &  \Bigg \{\lambda f(\cdot - a) + (1-\lambda) f(\cdot -b) , \  f \in \mathcal{SLC}, \  f(0) \in [1/M, M],  \lambda \in [0,1],  \\
  && \hspace{0.2cm}\ (a,b)  \in [a_0 - \delta, a_0+ \delta] \times [b_0 - \delta, b_0 + \delta] \Bigg \}.
\end{eqnarray*}
Here, the parameters $a_0$ and $b_0$ play the role of the true location
shifts $u^0_1$ and $u^0_2$. Consistency of the estimates
$\uestI$ and $\uestII$
ensures that they are stay within distance $2 \delta$ from the truth with
increasing probability. Also, uniform consistency of the log-concave MLE,
$\widehat f_n$, on continuity sets of $f^0$ implies consistency at the point
0 (the common mode of $f^0$ and $\widehat f_n$). Thus, we can find $M > 0$
such that $\widehat f_n(0) \in [1/M, M]$ with increasing probability. The
following proposition gives a bound on the bracketing entropy for the class
$\mathcal G$.

\medskip

\begin{prop}\label{Entropy}
  For $\epsilon \in (0, \epsilon_0]$, we have that
  \begin{eqnarray*}
    \log N_{[ \ ]}(\epsilon,\mathcal G,H)  \lesssim \frac{1}{\sqrt \epsilon}
  \end{eqnarray*}
  where $\epsilon_0$ and $\lesssim$ depend only on $a_0, b_0, \delta$ and $M$.
\end{prop}

\medskip

\bigskip

\par \noindent Now, we are ready to state our main theorem.  We find a rate
of convergence of at least $n^{-2/5}$ in the $L_1$ norm, both for $\widehat
f_n $ and $\widehat g_n$.  Although we consider $L_1$ and
\cite{ButuceaVand2014} consider $L_2$ distance, the rate of
Theorem~\ref{Rates} is an improvement over the corresponding $L_2$ rate $n^{-
  (2\beta -1 ) /(4\beta)}$ of \cite{ButuceaVand2014} whenever $\beta < 5/2$.
(Note that (log-)concave functions are Lebesgue-almost-everywhere twice
differentiable by Alexandrov's theorem \citep{niculescu2006convex}, so
roughly correspond to $\beta$ being $2$ or larger.)

\medskip

\begin{theorem}\label{Rates}
  Let $\widehat f_n$ and $\widehat g_n$ be again the MLE's of the zero-symmetric
  log-concave component and mixed density respectively. If $\sqrt n (\uestI -
  u^0_1) = O_p(1)$, $\sqrt n (\uestII - u^0_1) = O_p(1)$, and
  $\sqrt n (\piest - \pi^0) = O_p(1)$, then %
  \begin{eqnarray*}
    L_1\big(\widehat f_n, f^0\big) = O_p\big(n^{-2/5}\big), \ \textrm{and} \ \  L_1\big(\widehat g_n, g^0\big) = O_p\big(n^{-2/5}\big)
  \end{eqnarray*}
  where $L_1(d_1, d_2) = \int_{\RR} \vert d_1(x) - d_2(x) \vert dx$.

\end{theorem}

\bigskip

\begin{figure}
  \begin{minipage}[t]{0.5\textwidth}
    \includegraphics[width=\textwidth]{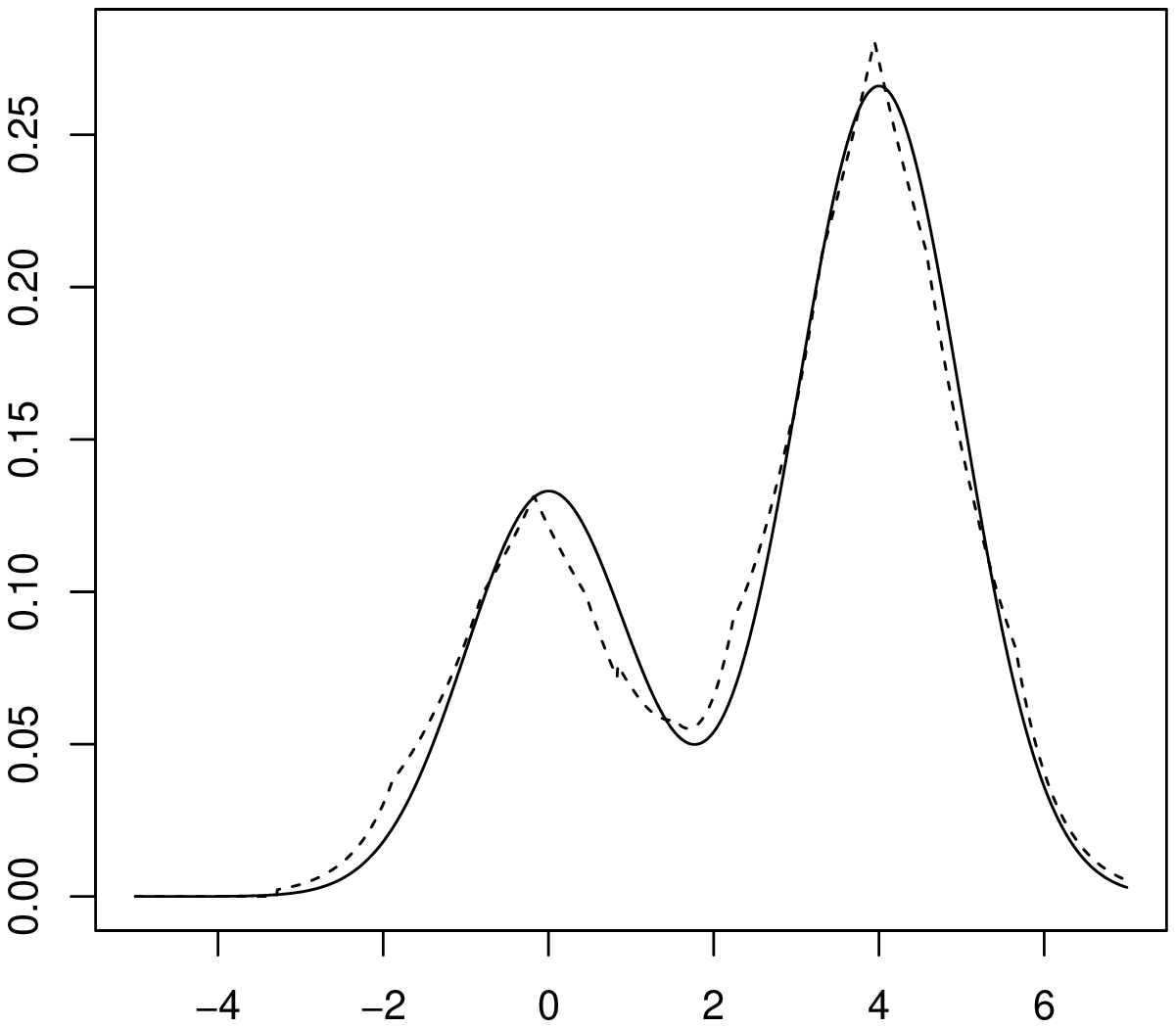}
  \end{minipage}
  \begin{minipage}[t]{0.5\textwidth}
    \includegraphics[width=\textwidth]{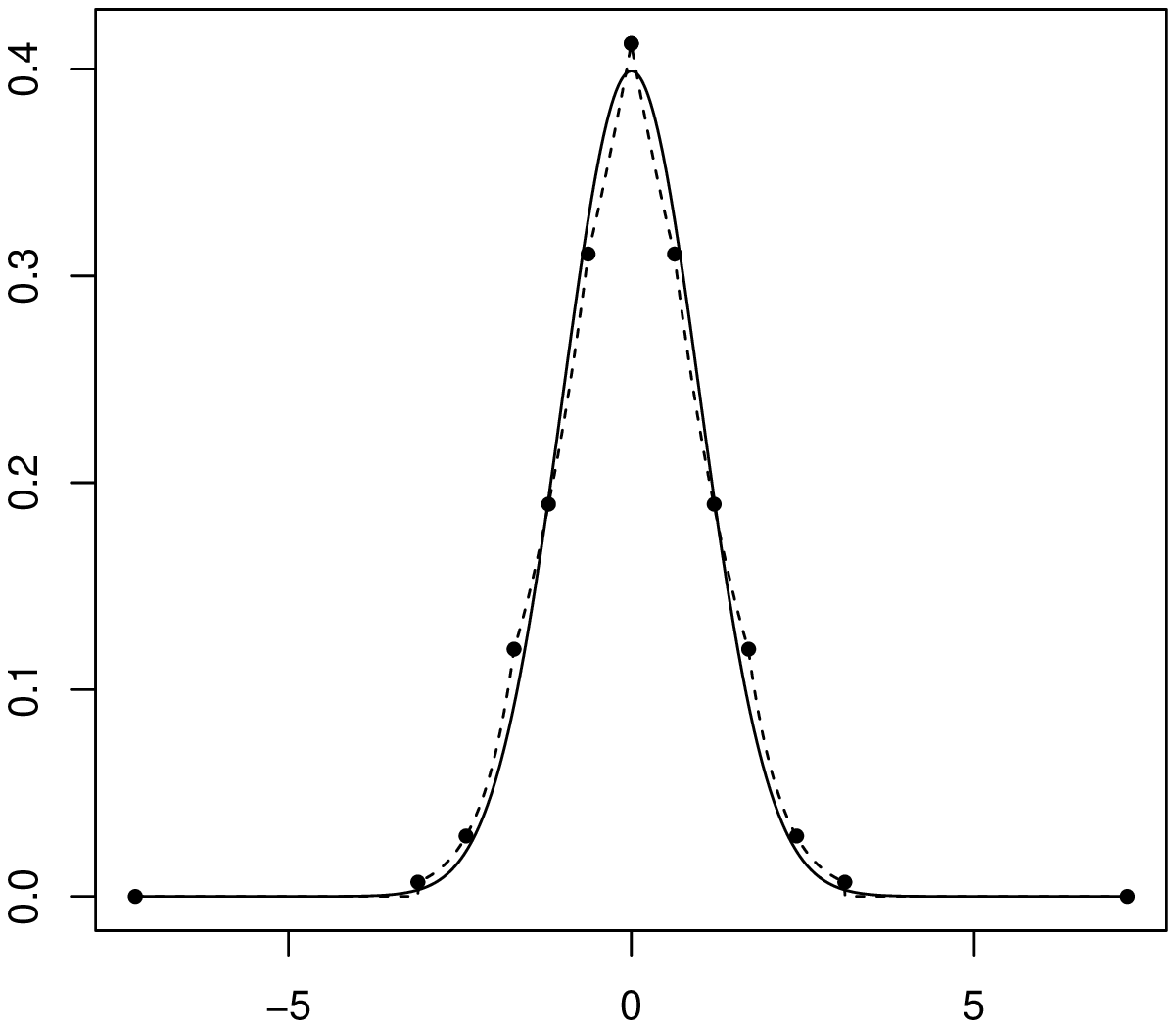}
  \end{minipage}
  \caption{Left: plot of the true density of $(1/3) \mathcal{N}(0,1) +(2/3) \mathcal{N}(4,1)$ (solid line) and its log-concave estimator $\widehat{g}_n$ (dotted line).  Right: plot of the density of $\mathcal{N}(0,1)$ and its zero-symmetric log-concave MLE $\widehat f_n$. The MLE was based on $n=500$ independent data drawn from the mixture density $(1/3) \mathcal{N}(0,1) + (2/3) \mathcal{N}(4,1)$.}
  \label{FigMLestimates}
\end{figure}

To illustrate the theory, a simulated example is given in Figure
\ref{FigMLestimates}. The true zero-symmetric component $f^0$ is taken to be the
density of a standard Gaussian with mixing probability $\pi^0= 1/3$ and shift
locations $u^0_1 = 0$ and $u^0_2 = 4$.  The plot on the left (right) shows
our MLE of $g^0$ ($ f^0 $) based on a sample of size $n=500$. The bullets in
the right plot depict the knot points of the zero-symmetric log-concave MLE, that
is, the points where the logarithm of the log-concave MLE changes its slope.

\section{Testing and clustering}\label{NumSection}

\cite{bordes2006}, \cite{hunteretal2007} and \cite{ButuceaVand2014} propose
three different ways of estimating the mixture parameters $\pi^0, u^0_1$ and
$u^0_2$. As we are interested here in $\sqrt n$-consistent estimators of
these parameters, we prefer the work by \cite{hunteretal2007} and
\cite{ButuceaVand2014}.
Also, due to some numerical instabilities encountered when computing the
estimators proposed by \cite{ButuceaVand2014}, we adopt the approach of
\cite{hunteretal2007} which has been already implemented in \proglang{R}; one
could either used the code posted at
\url{http://www.stat.psu.edu/~dhunter/code} or the function in \pkg{the
  mixtools} package. The latter option was kindly brought to the attention of
the first author by David Hunter in a private communication.

Once the estimates $\piest$, $\uestI$ and $\uestII$ of $\pi^0$, $u^0_1$ and $u^0_2$ are computed,
we maximize
\begin{eqnarray*}
  f^+ \mapsto   \frac{1}{n} \sum_{i=1}^n \log\left(\piest f^+(Z_{2i-1}) + (1-\piest) f^+(Z_{2i}) \right)
\end{eqnarray*}
where $f^+ = 2 f 1_{[0, \infty)}$, $Z_{2i-1} = \vert X_i - \uestI
\vert$ and $Z_{2i} = \vert X_i - \uestII \vert$ for $i = 1, \ldots,
n$.  This is equivalent to maximizing (\ref{LL}), which we have shown to admit
a maximizer.  Since the log-likelihood is not concave and it is not clear how
to maximize it directly, we will appeal to the EM algorithm \citep{MR0501537}.  Although we are
fixing the parameters $\pi, u_1$, and $u_2$ we
may still introduce the standard-in-mixture-models complete data of $(X_i,
\Delta_i)$, where $\Delta_i \sim $ Bernoulli$(\pi)$ and $X_i| \lb
\Delta_i=1 \rb \sim f^0(\cdot - \uestI)$ and $X_i | \lb \Delta_i=0
\rb \sim f^0(\cdot - \uestII)$.  An iteration of the EM algorithm in
this setup is then, given an estimate $\hat f^{+,(r)}_n $, to compute
\begin{equation}
  \label{eq:EM-argmax}
  \hat f^{+,(r+1)}_n :=
  \argmax_{f^+}
  \sum_{i=1}^n \left\{\hat p^{(r)}_{n,i} \log f^+(Z_{2i-1}) + (1-\hat
    p^{(r)}_{n,i}) \log f^+(Z_{2i}) \right\}
\end{equation}
where the argmax is over log-concave densities on $[0,\infty)$ with mode at
$0$, and
\begin{equation*}
  \widehat p^{(r)}_{n,i}  = \frac{\piest \widehat{f}^{+,(r)}
    (Z_{2i-1})}{ \piest \widehat{f}^{+, (r)} (Z_{2i-1})  + (1- \piest) \widehat{f}^{+,(r)} (Z_{2i})}
\end{equation*}
is the conditional expectation of $\Delta_i$ given $X_i$.  %
To initialize the EM algorithm, we start
with the density of a centered Gaussian distribution with variance equal to
the estimate given in formula (11) of \cite{hunteretal2007} for the true
variance of the zero-symmetric component, that is %
$\frac{1}{n} \sum_{i=1}^n (X_i - \bar{X}_n)^2 - \piest (1- \piest) (\uestII- \uestI)^2,$
or $1$ if this estimate is negative  (this may occur for moderate sample sizes).
The argmax in \eqref{eq:EM-argmax} can be computed by
the \proglang{R} package
\pkg{logcondens.mode}.

\subsection{Testing the absence of mixing}\label{subsec:test}

Recall that the mixing model we consider in this paper is given by
\begin{eqnarray*}
  \label{eq:MixtureModel-test}
  g^0 = \pi^0 f^0(\cdot - u^0_1) + (1-\pi^0) f^0(\cdot - u^0_2)
\end{eqnarray*}
with $f^0$ a log-concave zero-symmetric density on $\RR$, $\pi^0 \notin \{0, 1/2,
1 \}$ and $u^0_1 < u^0_2$.  We now use our log-concave MLE to test for the
absence of mixing, i.e.\ to test for the null hypothesis that $u_1^0=u_2^0$,
against the alternative that $u_1^0 \ne u_2^0$ and $\pi^0 \ne 1/2$ under the assumption that $f^0$ is zero-symmetric and log-concave.

To test for mixing we consider the likelihood ratio statistic. Under
the null hypothesis, we take the estimator of the true density to be equal to
the log-concave MLE which is symmetric around the median of the data. If
$\widehat g_n^0$ denotes this estimator, then our test statistic is given by
\begin{equation}
  \label{eq:LRS}
  \Lambda_n
  = \frac{ \prod_{i=1}^n \widehat g_n(X_i) }{\prod_{i=1}^n \widehat g_n^0(X_i)
  }.
\end{equation}
The null hypothesis is then rejected when $\Lambda_n$ is too large.  We use
the null hypothesis estimator to find critical values; that is, we
bootstrap from the symmetric log-concave estimator $\widehat g^0_n$.  The critical values of $\Lambda_n$ are then
computed in the usual way: based on the bootstrapped samples from $\widehat
g_n^0$, we compute the estimators of the mixing probability and mixture
locations and the corresponding MLE $\widehat g_n$. The order
statistics of the bootstrapped values of the likelihood ratio are then
obtained to compute upper empirical quantiles of a given order.
We also compare our test for mixing (hereafter referred to as the LR test) to the following procedures:

\begin{itemize}

\item the naive symmetric bootstrap (NSBS):  we re-sample with replacement $n$ random variables $Z^*_1, \ldots, Z^*_n$ from $\{ \pm \vert X_1 - \widehat m_n \vert, \ldots, \pm  \vert X_n - \widehat m_n \vert\}$ with  $\widehat m_n$ the median of $X_1, \ldots, X_n$ and set $X^*_i = \widehat m_n + Z^*_i, i = 1, \ldots, n$.  Then, the bootstrapped estimators of the location mixture of \cite{hunteretal2007}, $u^*_1$ and $u^*_2$, are computed based on $Y^*_i, i=1, \ldots, n$. We repeat this procedure $B$ times and compute the empirical $(1-\alpha)$-quantile of the distribution of $u^*_2 - u^*_1$.  The null hypothesis is rejected if the observed $\uestII - \uestI$ is larger than this quantile.

\item the naive symmetric bootstrap based on symmetric kernel density estimation (NSBSKDE): the method is similar to the one described above except that a standard kernel density estimator is fitted to $\widehat m_n \pm \vert X_i - \widehat m_n \vert$ and $X^*_1,  \ldots, X^*_n$ are now drawn from the fitted estimator at each bootstrap iteration.

\item the likelihood ratio based on symmetric kernel density estimation  (LRSKDE):  two
  kernel density estimators are computed, one under the full model, that is, based
  on $X_1, \ldots, X_n$, and one under the null model, that is, based on
  $\widehat m_n \pm \vert X_i - \widehat m_n \vert$. The likelihood ratio of
  these estimators is then computed. Bootstrap samples are obtained by
  simulation from the kernel estimator under the null hypothesis and then the
  empirical $(1-\alpha)$-quantile of the likelihood ratio is thereby computed.  The null hypothesis is rejected if the observed likelihood ratio is larger that this quantile.
\end{itemize}

Note that the NSBS provides a comparison procedure not based on density
estimation of the components.
In assessing the power, we take the true zero-symmetric component $f^0$ to be
one of the following distributions:
(1)  a standard Gaussian,
(2)  a double exponential, and
(3)  a uniform on $[-1,1]$,
Also, we take the true parameters to be
$\pi^0 \in \{0.20,  0.40 \}$ and
$(u^0_1, u^0_2) \in \{ (0,0), (0, 1), (0,3) \}$.
We give the estimated probability of rejecting the null hypothesis based on
$R=500$ replications with $B=49$ bootstrap samples in
Table~\ref{PowerLRTRn250}, for $n=250$.
The simulation results show the LR and LRSKDE tests are both outperforming the NSBS and NSBSKDE with power nearly equal or equal to $1$ for the well-separated mixtures. However, all the considered tests  seem to have a level larger than the specified level $\alpha=0.1$ for the uniform distribution.  Further simulations, which we do not report here, show that this improves when the sample size is increased to $n=500$. Note that the mixtures with mixture probability $\pi^0 =.4$ are more difficult to distinguish than those with $\pi^0 = .2$.  This is to be expected as the former mixtures are close to being symmetric around the mid-point $(u^0_1 + u^0_2)/2$.

It would be interesting to know whether the  level of our testing procedure based on the bootstrapped likelihood ratio test equals the theoretical level. The problem is however far from being trivial. Deriving the asymptotic level for example would require establishing the limit distribution of our statistic under the null hypothesis and also showing that it admits a continuous cumulative distribution function. Establishing such results requires a thorough study of the global asymptotics of the log-concave MLE. As this is outside the scope of this paper, the question remains open.

\begin{table}
  \centering
  \small %
  \begin{tabular}{ccl||ccccc}
    \hline
    \raisebox{-1.3ex}{\bf  Distribution} &  \raisebox{-1.3ex}{$\boldsymbol{\pi^0}$}  &  \backslashbox{{\bf Test}}{$\boldsymbol{u^0_2-u^0_1=}$}  & \raisebox{1.3ex}{$ 0$} & \raisebox{1.3ex}{$1$} & \raisebox{1.3ex}{$3$} \tabularnewline
    \hline
    \hline
    \multirow{8}{*}{$\mathcal{N}(0,1)$}
    & \multirow{4}{*}{$ .2$} &  LR  &  $0.11$ & $0.15$  & $1.00$  \tabularnewline
    \cline{3-6}
    & & NSBS  & $0.06$  &  $0.07$  &   $0.34$  \tabularnewline
    \cline{3-6}
    & & NSBSKDE &  $0.07$  &  $0.08$  &  $0.40$     \tabularnewline
    \cline{3-6}
    & & LRSKDE  &  $0.11$  & $0.11$   &   $1.00$   \tabularnewline
    \cmidrule[1pt]{2-6} %
    & \multirow{4}{*}{$ .4$} &  LR  & *  &  $0.14$  & $0.86$  \tabularnewline
    \cline{3-6}
    & & NSBS  & *  & $0.05$  &  $0.28$  \tabularnewline
    \cline{3-6}
    & & NSBSKDE  & *  &  $0.08$  &  $0.26$  \tabularnewline
    \cline{3-6}
    & & LRSKDE  & *  & $0.11$  & $0.87$   \tabularnewline
    \cmidrule[1pt]{1-6} %
    \multirow{8}{*}{$\mathcal{L}(1)$}
    & \multirow{4}{*}{$ .2$} &  LR  &  $0.11$  & $0.22$  &  $0.99$ \tabularnewline
    \cline{3-6}
    & & NSBS  & $0.11$  &   $0.01$  &   $0.01$  \tabularnewline
    \cline{3-6}
    & & NSBSKDE  & $0.13$  &  $0.11$   &   $0.02$    \tabularnewline
    \cline{3-6}
    & & LRSKDE  & $0.12$  &  $0.19$   &  $0.98$  \tabularnewline
    \cmidrule[1pt]{2-6} %
    & \multirow{4}{*}{$ .4$} &  LR  & *  &  $0.14$  &  $0.89$  \tabularnewline
    \cline{3-6}
    & & NSBS  & *  &  $0.08$  &  $0.12$   \tabularnewline
    \cline{3-6}
    & & NSBSKDE & *  &  $0.08$  &   $0.03$    \tabularnewline
    \cline{3-6}
    & & LRSKDE  & *  &  $0.18$   &  $0.86$   \tabularnewline
    \cmidrule[1pt]{1-6} %
    \multirow{8}{*}{$\mathcal{U}[-1,1]$}
    & \multirow{4}{*}{$ .2$} &  LR  & $0.19$  &  $0.92$  &   $1.00$ \tabularnewline
    \cline{3-6}
    & & NSBS  &  $0.11$ &  $0.07$  &  $0.16$   \tabularnewline
    \cline{3-6}
    & &NSBSKDE  &  $0.11$ &  $0.08$  &  $0.14$     \tabularnewline
    \cline{3-6}
    & & LRSKDE  & $0.25$  & $0.99$   &  $1.00$  \tabularnewline
    \cmidrule[1pt]{2-6} %
    & \multirow{4}{*}{$ .4$} &  LR  & *  & $0.60$  &  $1.00$  \tabularnewline
    \cline{3-6}
    & & NSBS  & *  &  $0.09$  &  $0.78$   \tabularnewline
    \cline{3-6}
    & & NSBSKDE  & *  &  $0.07$  &  $0.75$     \tabularnewline
    \cline{3-6}
    & & LRSKDE & *  &  $0.54$  &  $1.00$   \tabularnewline
    \cmidrule[1pt]{1-6} %
  \end{tabular}
  \caption{Values of the bootstrapped power for LR, NSBS, NSBSKDE and LRSKDE
    tests when the true density is $\pi^0 f^0(\cdot - u^0_1) + (1-\pi^0) f^0(\cdot - u^0_2)$,
    where $u^0_1=0$, $u^0_2-u^0_1 \in \{0,1,3\}$, $\pi ^0 \in \{0.2, 0.4\}$, and $f^0$
    is one of the zero-symmetric log-concave densities  shown in the first column. The nominal level is $\alpha=0.1$. The  sample  size is
    $n=250$,  the number of bootstraps and the number of replications were taken to be $B=49$ and $R=500$, respectively.  The common value of the power under $H_0$ is replaced by  \lq\lq *\rq\rq.}
  \label{PowerLRTRn250}
\end{table}

\begin{table}
  \small %
  \centering
  \begin{tabular}{c||cccc}
    & G & HG & SLC &KDE \tabularnewline \hline
    $\mathcal{N}(0,1)$ & $163 \, (0.44) $ & $170 \, (0.30) $ & $170 \, (0.30)$ & 182 (0.88) \tabularnewline
    $\mathcal{L}(1)$ & $223 \, (0.64) $  & $173 \, (0.32) $  & $174 \, (0.28) $ & 154 (0.82) \tabularnewline
    $\mathcal{U}(-1,1)$ & $144 \, (0.15) $  & $66  \, (0.11) $ &  $55 \, (0.10) $ & 105 (0.28) \tabularnewline
    \hline
  \end{tabular}
  \caption{
    Comparison of the four different clustering methods given by the column labels, see text for more
    details. The reported numbers are the average number of
    misclassifications out of $n=500$ samples over $R=5000$ replications
    under each of the three log-concave densities in the left column.  Here  $u^0_2-u^0_1=1$ and $\pi^0 = .2$. The numbers in parentheses are the corresponding standard errors. %
  }
  \label{tab:class}
\end{table}

\subsection{Gaussian versus symmetric log-concave clustering}

We now consider the problem of clustering, i.e., of assigning to each
observation in a dataset a label without being given any ``training''
labels.  We will assume that the data can be clustered into two groups,
which we will do by by fitting the two-component mixture
\eqref{eq:MixtureModel} and assigning a label to an observation $X$ based on
whether our estimate of the posterior probability
\begin{eqnarray}
  \label{eq:posteriorprob}
  \frac{\pi^0  f^0(X - u^0_1)}{\pi^0  f^0(X - u^0_1) + (1-\pi^0)  f(X - u^0_2) }
\end{eqnarray}
is greater than $1/2$ or not.

We fit the mixture three different ways.  In the first basic approach, labeled ``G'', we maximize the
likelihood \eqref{LL} under the assumption that the component $f$
is a normal density.  We use the EM algorithm  to maximize
the likelihood.  Our next two approaches both use the method of
\cite{hunteretal2007} to estimate the mixture components $u^0_1,u^0_2$, and
$\pi^0$.  Then we either fit the components using a {\em G}aussian density (denoted
``HG''), with variance estimate also given by \cite{hunteretal2007}, or we
use the symmetric log-concave density estimator (denoted ``SLC'') for the
components.  The fourth approach is based on the estimators of \cite{hunteretal2007} and the kernel density estimator based on the inversion formula given in (9) by \cite{bordes2006} where we truncate the infinite sum at some large integer $K > 0$. Precisely, let $\bar{g}_n$ be a standard density estimator of the mixed density $g^0$. Then, the KDE of $f^0$ we use is given by
$\max(1/2 (\bar{f}_n(x) + \bar{f}(-x)), 0)$ where
\begin{eqnarray*}
  \bar{f}_n(x) = \sum_{k=0}^K \left( \frac{-\piest}{1-\piest}  \right)^k  \bar{g}_n\big(x  + \uestII  + k (\uestII - \uestI)\big).
\end{eqnarray*}
We should note that this formula is only valid, when $\piest < 1/2$.  Hence,  $1-\piest$ and $(\uestII, \uestI)$ should replace $\piest$ and $(\uestI, \uestII)$ When $\piest \ge 1/2$.

\smallskip

We record the average  missclassification count when the true density
is one of the densities in the left column of Table \ref{tab:class}. In all cases the number of replications is $R=5000$,
the sample size is $n=500$, $u^0_2-u^0_1=1$ and $\pi^0=.2$.

The performances of the four approaches were then compared, and the
results are reported in Table~\ref{tab:class}.  The KDE approach does clearly worse than the three other methods. The SLC
outperforms HG by $16$\% when the true density is $\mathcal{U}[-1,1]$.  In the other cases, they perform similarly.  All four methods define the two cluster regions by dividing the real line into two
half-lines.  The HG and SLC methods have the same mixture components so the shape of
the component density estimates have to be dramatically different (e.g.,
uniform instead of normal) in order to noticeably change the results;
note this somewhat deceiving outcome is not totally in contradiction with
the finding of \cite{cule_10} about the performance of their
two-dimensional log-concave classifier applied to the Breast cancer data of
Wisconsin; see \cite{cule_10} for details. The authors found that the
log-concave MLE reduces the percentage of misclassification from 10.36\%
obtained for the Gaussian estimate to only 8.43\% for that particular data
set. The posterior probabilities of cluster membership, which can
be used as a measurement of uncertainty, can also differ noticeably between
the HG method and our SLC method.

\section{Data application}
\label{sec:data}

In this section, we apply our new estimation approach to two different datasets.

\subsection{Old Faithful data}
The data to which we first apply our estimation procedure are the times, in
minutes, between eruptions of the Old Faithful geyser in Yellowstone National
park.  There are many forms of the Old Faithful data.  As far as we know, the
oldest version of the data was collected by S.\ Weisberg from R.\ Hutchinson
in August 1978. %
The data we analyze were collected between August 1 and August 15, 1985
continuously, and are from \cite{Azzalini:1990ep}.  The following explanation
from \cite{Weisberg:2005ja} motivates interest in the data:
\begin{quote}
  Old Faithful Geyser is an important tourist attraction, with up to several
  thousand people watching it erupt on pleasant summer days.  The park service
  uses data like these to obtain a prediction equation for the time to the next
  eruption.
\end{quote}

\vspace{-0.1cm}
\begin{figure}[h!]
  \begin{center}
    \includegraphics[width=1.1\textwidth]{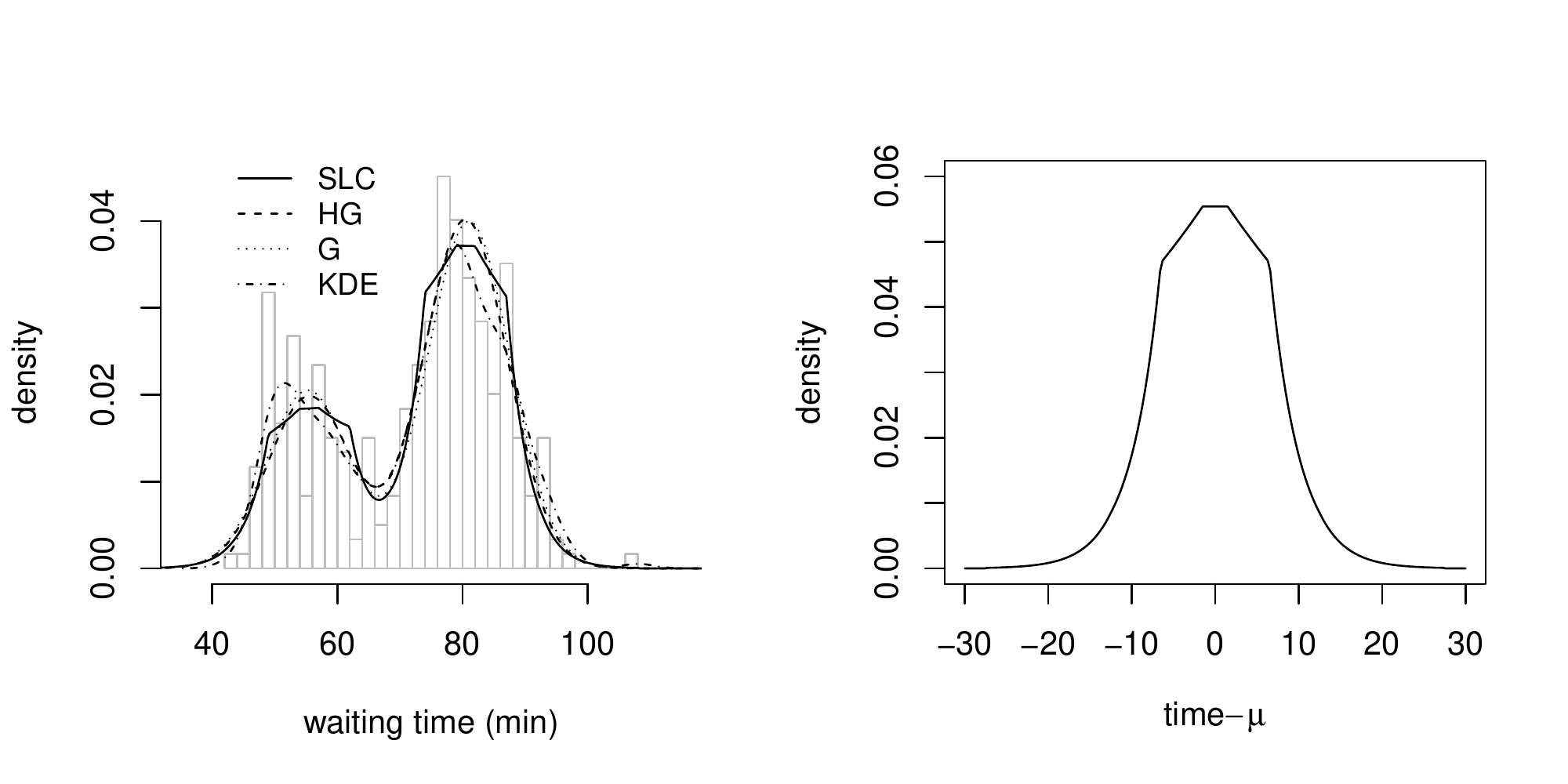}
  \end{center}
  \caption{Time between eruptions of Old Faithful Geyser (min).  The ``SLC''
    and ``HG'' estimates both use the method of \cite{hunteretal2007} to
    estimate the mixture parameters $u^0_1, u^0_2,$ and $\pi^0$.  ``SLC'' then fits with symmetric
    log-concave components and ``HG'' fits with Gaussian components.  The ``G''
    estimate is the maximum likelihood estimate of a mixture model with two
    Gaussian components with equal variances.  The "KDE" is a standard kernel density estimator with an optimal bandwidth.}
  \label{fig:geyser}
\end{figure}

In Figure~\ref{fig:geyser}, we have two plots related to the Old Faithful
data. The plot on the left depicts a descriptive histogram of the
data with around $30$ bins (which is too many for optimal estimation)
along with the plots of four mixture density estimates. The ``SLC'' ({\em s}ymmetric {\em l}og-{\em
  c}oncave) estimate is the mixture model where $u^0_1,u^0_2$ and $\pi^0$ are
estimated using the method of \cite{hunteretal2007}, and then the components
are estimated using our symmetric log-concave estimator.  The ``HG'' ({\em
  H}unter et al.\ and {\em G}aussian components) estimate is given by again
using the method of \cite{hunteretal2007} to find estimates of the mixture parameters
whereas the nonparametric components are taken to be Gaussian components (the same Gaussian density for both
mixture components).  The estimates for $u^0_1 ,u^0_2,$ and $\pi^0$ given by
\cite{hunteretal2007} are $55.5$, $80.5$, and $.33$, respectively. The ``G''
({\em G}aussian) estimate in the plot is based on simply using a Gaussian
mixture model with two components with equal variances.  Assuming equal
variances forces the two components to be identical, which makes the model
analogous to the others. In this case, we estimated $u^0_1, u^0_2, $ and
$\pi^0$ by the EM algorithm \citep{MR0501537}, with estimated values of
$55.3$, $81.0$, and $.339$. The normal components are slightly more peaked
than the log-concave ones, but the overall fit is fairly similar; in large
part this is because the locations and weights are very similar. Finally, the \lq KDE\rq \ is a standard kernel density estimator with an optimal bandwidth.

The plot on the right is that of the zero-symmetric log-concave component,
centered at $0$, used in the mixture density.  As expected from the known
theoretical properties of this estimator, it has a flat interval about the origin,
and is the exponential of a concave piecewise linear and zero-symmetric function.

\subsection{Height data}

We next examine $1766$ human height observations.  We look at the
heights of the population of Campora, a village in the south of Italy.  This
population is studied by the ``Genetic Park of Cilento and Vallo di Dano
Project'' \citep{geneticpark}, which is interested in identifying
geographically and genetically isolated populations. Such populations are of
particular interest because in addition to ``genetic homogeneity,'' they have
a ``uniformity of diet, life style and environment.''  These homogeneities
are valuable in the study of genetic risk factors for complex pathologies
such as ``hypertension, diabetes, obesity, cancer, and neurodegenerative
diseases,'' by allowing for a ``simplification of the complexity of genetic
models'' involved, because  %
of the population's homogeneity \citep{geneticpark}.

\cite{colonna2007campora} provide evidence that this population is indeed
genetically isolated.  Because of this feature, the distribution of heights
of this population is not necessarily the same as that of the global
population at large, so estimating its distribution is of interest.  Height
data are often modeled as mixtures of two components, corresponding to the
two sexes.

\begin{figure}[h!]
  \begin{center}
    \includegraphics[width=1\textwidth]{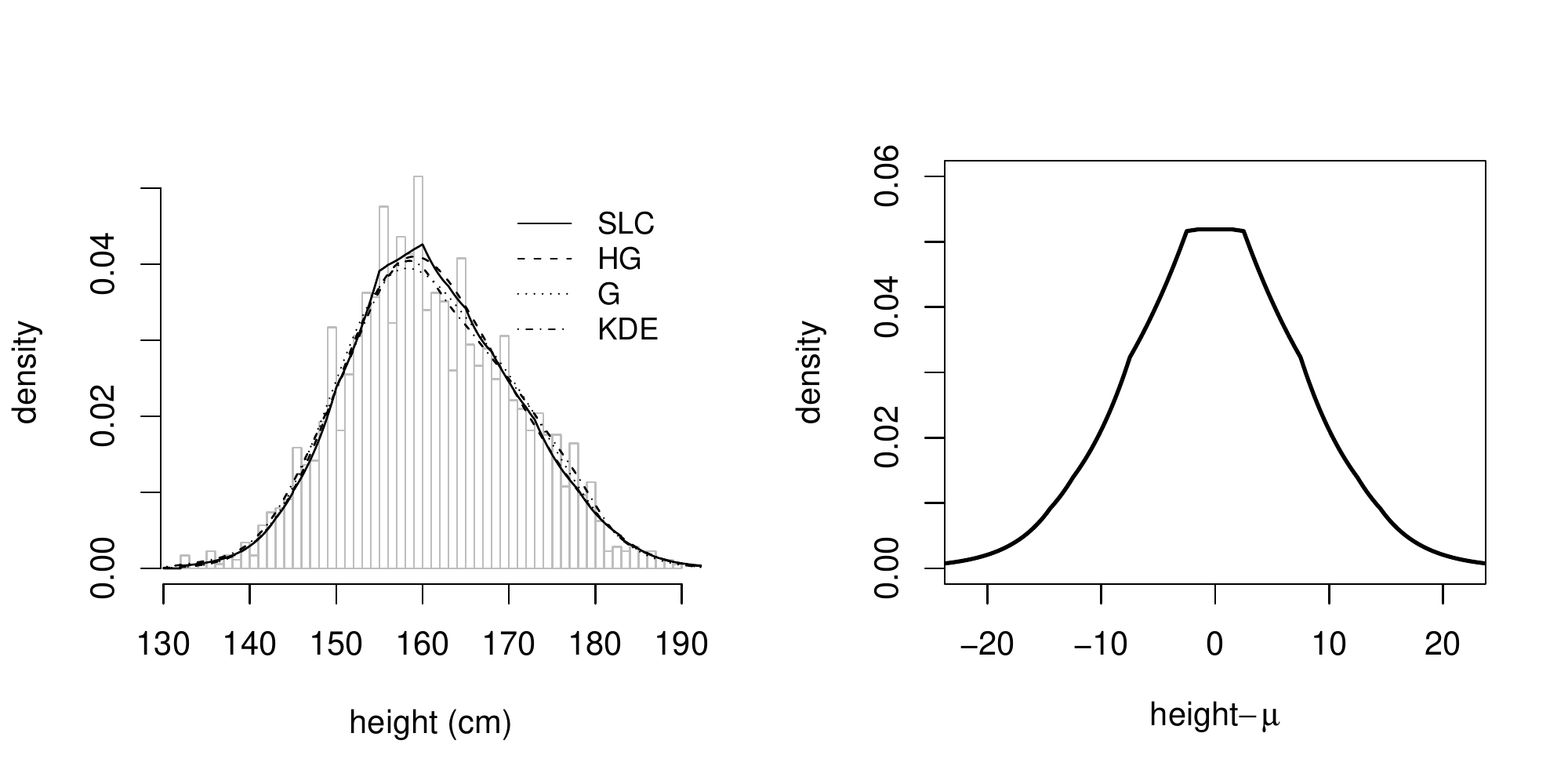}
  \end{center}
  \caption{ Height data of the population of Campora . The ``SLC''
    and ``HG'' estimates both use the method of \cite{hunteretal2007} to
    estimate the mixture parameters $u^0_1, u^0_2,$ and $\pi^0$.  ``SLC'' then fits with symmetric
    log-concave components and ``HG'' fits with Gaussian components.  The ``G''
    estimate is the maximum likelihood estimate of a mixture model with two
    Gaussian components with equal variances.  The \lq KDE\rq \ is a standard kernel density estimator with an optimal bandwidth.}
  \label{fig:height}
\end{figure}

We present plots related to the height data in Figure~\ref{fig:height}.
The height data do not exhibit multi-modality, but
two-component mixtures still fit the data well.  The three approaches that
we consider fit
similarly, but the log-concave components are able to capture a bit more
asymmetry near to the mode.

\begin{figure}[h!]
  \begin{center}
    \includegraphics[width=1\textwidth]{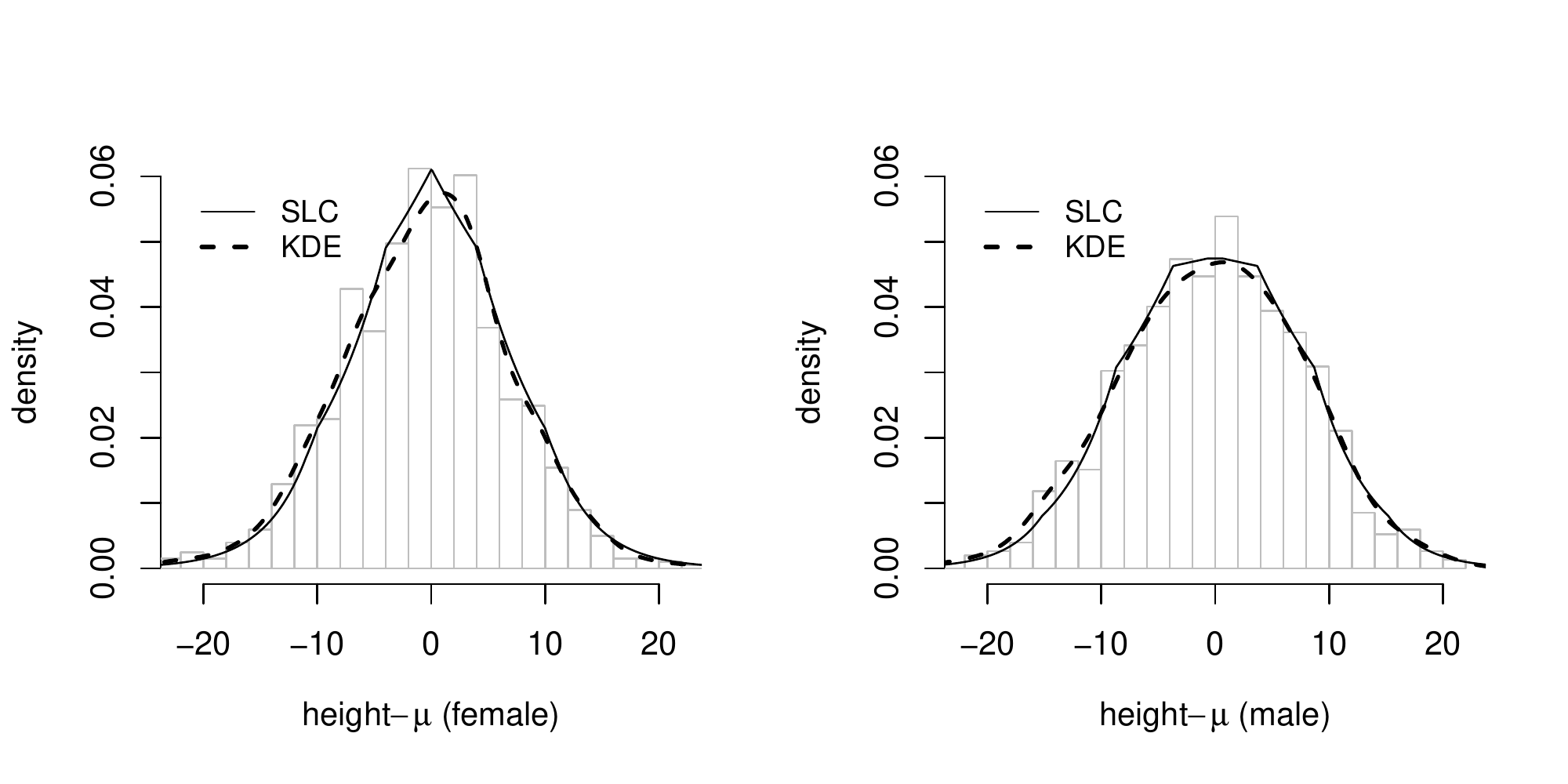}
  \end{center}
  \caption{Descriptive histograms of the height data for women (left) and men (right) after centering around the median. The ``SLC" is the log-concave MLE of the true density constrained to have mode at $0$. The ``KDE" is a standard kernel density estimator with an optimal bandwidth.}
  \label{fig:height_split_components}
\end{figure}

The plot on the right includes the mixture component density (labeled
``All''), in black.  The data include the sex of each individual, so, using
this extra information we can also estimate the true component densities
separately: the zero-symmetric log-concave
density estimate can be compared to the estimates of the density of the
heights for either sex considered alone.
Figure \ref{fig:height_split_components} shows the plots of the
(descriptive) histograms
of the heights for men and women and fitted standard kernel density
estimators.
The assumption of symmetry of the distribution of the heights for each of the
genders seems to be reasonable to make.  The observed proportion of women is
$ 0.57$, whereas the observed medians of the heights for  women and men were
found to be $156.0$ and $168.7$,  respectively.
Using the estimation method of \cite{hunteretal2007}, we found  $\piest=0.72,
\uestI = 157.5,$ and $\uestII=170.5$.   Here $\piest$ and $\uestI$ correspond
to the component for women.
The components estimated by using the labels for men and women differ from
that using the mixture model without the labels especially towards the
center.  We believe that this is essentially due to the difference between the estimates of the mixture parameters $\pi^0$, $\mu_1^0$ and
$\mu_2^0$ obtained by ignoring or using the information available about the gender. In the latter case, the locations are estimated by the respective medians.  It does appear that the
distributions of heights of men and women are somewhat different near those
centers, with women having a more peaked density and
men having a flatter one.  Thus, in the mixture model, without using the
labels, the component density estimate is somewhere in between the two
shapes.

\section{Conclusions}

The goal in this paper is to make use of the log-concavity constraint to
estimate the unknown density component in a semi-parametric location mixture
model assuming that this unknown density is symmetric around the origin. The
first motivation for choosing this approach is that many densities are log-concave. The second one is to build an estimation
procedure that does not depend on a tuning parameter.  Our log-concave MLE is computed by maximizing the log-likelihood function after estimating the mixture parameters using the approach of \cite{hunteretal2007}.  The computation is  easily implementable using the EM algorithm in combination
with an active set algorithm already implemented in the \proglang{R} package
\pkg{logcondens.mode}.

As already mentioned, our method is not advocated for heavy-tailed densities.  In
such cases, other shape constraints may be more appropriate,
specifically, $s$-concavity, as studied in \cite{mizera_09} and
\cite{dosswellner13}. Unfortunately, the  theory of estimators of
$s$-concave densities is less developed than that of log-concave
MLEs, which remains a barrier to using $s$-concavity in our current
context.  %

Finally,  \cite{hunteretal2007} give
sufficient conditions on the mixing probabilities and mixture locations for
the model to be $3$-identifiable.
In this case, the mixture parameters can still be computed using
the method of \cite{hunteretal2007}, and the log-concave MLE can be computed
as described in this paper.  However, it is not immediate in that case whether the
same proof approaches would still yield the same rate of convergence.
Recently, \cite{BalabdaouiButucea2014}
proved that the number of components $k$, the mixture parameters and the
unknown density are identifiable provided that the density is P\'olya
frequency (of infinite order) such that its expectation is equal to zero.
For a precise definition of P\'olya frequence functions, we refer to
\cite{schoenberg51}. The obtained identifiability result can be used of
course in the case of symmetry but it is certainly not a requirement. Note
that imposing the log-concave constraint in this setting is natural since
the class of P\'olya frequency functions is a subset of the log-concave
class as shown by \cite{schoenberg51}. One may argue that non-parametric
classes such as symmetric densities or P\'olya frequency functions with
expectation equal to zero are not large enough. However, it seems that
identifiability is hard to obtain if one allows
for large classes.

\appendix

\begin{supplement-file}
  \section{Supplement}
\end{supplement-file}

\begin{all-in-one-file}
  \section{Appendix}
\end{all-in-one-file}

\subsection{Main Proofs}

\par \noindent \textbf{Proof of Proposition \ref{Existence}.} \ To show
existence, we first start by proving that $\widehat \psi_n$ is necessarily a
piecewise linear function that is flat in the neighborhood of zero. Let $\psi
\in \mathcal{SC}$. Also, let $Z_{(1)}$ denote the first order statistic of
the transformed data $Z_i$'s defined in Proposition 2.1. Consider $\bar \psi$ the unique
concave function such that $\bar{\psi}(Z_i) = \psi(Z_i ), i=1, \ldots, 2n$,
$\bar{\psi}(0) = \bar{\psi}( Z_{(1)})$, $\bar \psi$ is piecewise linear
between the points $Z_i, i =1, \ldots, 2n$, $\bar{\psi}(t) = -\infty$ for $t
> Z_{(2n)}$ and $\psi(t) = \psi(-t)$ for $t \in \RR$. Clearly, $\bar \psi \in
\mathcal{SC}$ and admits the properties $\piest e^{\bar {\psi}(X_i - \uestI)} + (1-
\piest) e^{\bar{\psi}(X_i - \uestII)} = \piest e^{\psi(X_i - \uestI)} + (1- \piest)
e^{\psi(X_i - \uestII)}$ for $i=1, \ldots, n$ and $\bar{\psi}\le \psi$. This
implies that $\Phi_n(\bar{\psi}) \ge \Phi_n(\psi)$ and the logarithm of the
MLE has to be necessarily piecewise linear between the transformed data
points $Z_i$ and $-Z_i, i =1, \ldots, 2n$ with a flat piece in the
neighborhood of zero and support $[-Z_{(2n)}, Z_{(2n)}]$. Further, if a
maximizer $\widehat f_n = \exp(\widehat \psi_n)$ exists then it has to be a
density. This follows from the fact that $\widehat \psi_n + \epsilon \in
\mathcal{SC}$ for all $\epsilon \in \RR$, and
\begin{eqnarray*}
  0 = \lim_{\epsilon \to 0} \frac{\Phi_n(\widehat \psi_n + \epsilon) - \Phi_n(\widehat \psi_n)}{\epsilon}  = 1 - \int_{\RR} \exp(\widehat \psi_n).
\end{eqnarray*}
Given the results obtained above, a function $\psi \in \mathcal{SC}$ can be identified with the vector
\begin{eqnarray*}
  \underline{\psi} \equiv \Big(\psi(X_i - \uestI), \psi(X_i - \uestII)\Big)_{i=1}^n = \Big(\psi(Z_{2i-1}), \psi(Z_{2i})\Big)_{i=1}^n.
\end{eqnarray*}
Let $\mathcal{SC}_n$ denote the set of such vectors. In this second part of the proof, we will show that the maximization problem admits a solution. It is clear that
$\Phi_n$ is continuous on $\mathcal{SC}_n$. In the following, we will show that $\Phi_n$ is necessarily maximized  on a compact set. To this aim, consider $\underline{\psi}^{(p)}$ to be a
maximizing sequence, such that $\int \exp(\psi^{(p)}(t))dt =1$, and as $p \to
\infty$
\begin{eqnarray*}
  \max_{1 \le i \le n} \big \vert \psi^{(p)}(Z_{(i)}) \big \vert  \to \infty,   \ \textrm{and} \ \big \vert \psi^{(p)}(Z_{(i)}) \big \vert \to l_i \in [-\infty, \infty], \ \textrm{for $1\le i \le n$}
\end{eqnarray*}
where $Z_{(i)}$ denotes the $i$-th order statistic of $Z_1, \ldots, Z_{2n}$.

\bigskip

For $1 \le i \le n$ write $\psi^{(p)}_i := \psi^{(p)}(Z_{(i)})$. Suppose that  for all $i \in \{1, \ldots, n\}$
\begin{eqnarray*}
  \lim_{p \to \infty} \psi^{(p)}_i < \infty.
\end{eqnarray*}
This implies that there exists $j \in \{1, \ldots, n \}$ such that $ \lim_{p \to \infty} \psi^{(p)}_j  = -\infty$.   If there exists $i \in \{1, \ldots, n \}$ such that $Z_{(j)} = \max(\vert X_i - \uestI \vert, \vert X_i - \uestII \vert)$ then
\begin{eqnarray*}
  \lim_{p \to \infty} \log\left(\piest \exp(\psi^{(p)}(\vert X_i - \uestI \vert)) + (1-\piest)\exp(\psi^{(p)}(\vert X_i - \uestII \vert))\right) &\le & \lim_{p \to \infty} \psi^{(p)}_j
  \\
  & =& -\infty
\end{eqnarray*}
implying that $\lim_{p \to \infty} \Phi_n(\psi^{(p)}) = -\infty$ which is  in contradiction with the definition of $\psi^{(p)}$.

Let us assume now that for any $j$ such that $\lim_{p\to \infty} \psi^{(p)}_j
= -\infty$, there exists $i \in \{1, \ldots, n \}$ such that
\begin{equation}
  \label{Assump2}
  Z_{(j)} = \min(\vert X_i - \uestI \vert, \vert X_i - \uestII \vert) \ \textrm{and}
  \ \lim_{n \to \infty} \psi^{(p)}(\max(\vert X_i - \uestI \vert, \vert X_i -
  \uestII \vert) > -\infty.
\end{equation}
Suppose that $j=2n$ is the only integer for which this divergence is occurring. Let $i \in \{1, \ldots, n\}$ such that we have $Z_{(2n)} = \min(\vert X_i - \uestI \vert, \vert X_i - \uestII \vert)$. Without loss of generality,  assume that $\vert X_i - \uestI \vert < \vert X_i - \uestII \vert$. Put
\begin{eqnarray*}
  a^{(p)} = (1-\piest)\exp\left(\psi^{(p)}(\vert X_i - \uestII\vert)\right), \ \textrm{and} \ b^{(p)} = \psi^{(p)}_{2n-1}
\end{eqnarray*}
and call $a$ and $b$ their respective limits as $p \to \infty$. Note that both $a$ and $b$ are in $\RR$ by assumption.  Consider now the function defined by
\begin{eqnarray*}
  h(x)  & = & \log\left( \piest e^x + a \right) - 2 (Z_{(2n)} - Z_{(2n-1)}) \frac{\exp(b) - \exp(x)}{b - x}
\end{eqnarray*}
on $(-\infty, 0]$. Put $\Delta_{2n-1} = Z_{(2n)} - Z_{(2n-1)}$. We have that
\begin{eqnarray*}
  h'(x) = \frac{\piest \exp(x)}{\piest \exp(x) + a} -  \frac{2\Delta_{2n-1}}{b -x} \left(\frac{\exp(b) - \exp(x)}{b -x} - \exp(x)  \right)
\end{eqnarray*}
which has the same sign as
\begin{eqnarray*}
  &&\piest-  2\Delta_{2n-1} \exp(b) \Big(\piest \exp(x) + a\Big) \frac{1}{b -x} \left(\frac{\exp(-x) -\exp(-b)}{b -x} - \exp(-b)  \right) \\
  && =\piest - 2\Delta_{2n-1} \exp(b) \Big(\piest \exp(x) + a\Big) \frac{\exp(-x) -\exp(-b)}{(b -x)^2 }\left(1 - \frac{(b -x)\exp(-b) }{\exp(-x) - \exp(-b)} \right) \\
  && \to -\infty, \ \textrm{as $x \to -\infty$}.
\end{eqnarray*}
This implies that  as $p \to \infty$
\begin{eqnarray*}
  \log\left(\piest \exp(\psi^{(p)}(\vert X_i - \uestI \vert) + (1- \piest) \exp(\psi^{(p)}(\vert X_i - \uestII \vert))\right) - \int_{Z_{(2n-1)}}^{Z_{(2n)}} \exp(\psi^{(p)})(t) dt
\end{eqnarray*}
must be decreasing, and hence bigger values of the log-likelihood are to be obtained if $\psi^{(p)}(Z_{(2n)})$ were not divergent. This implies that $\psi^{(p)}_{2n}$ cannot diverge to $-\infty$. The same reasoning can be applied if (\ref{Assump2}) is satisfied by other integers $k \in \{1, \ldots, 2n-1\}$, and let $j$ be the smallest one. Let $i$ be such that $Z_{(j)} = \min(\vert X_i - \uestI \vert, \vert X_i -\uestII\vert )$. Note that $j > 1$. Also, note that by definition of $j$  $\lim_{p \to \infty} \psi^{(p)}(Z_{(j-1)}) > -\infty$.  Thus, we obtain the same conclusion as before, that is,  the term
\begin{eqnarray*}
  \log\left(\piest \exp(\psi^{(p)}(\vert X_i - \uestI \vert) + (1- \piest) \exp(\psi^{(p)}(\vert X_i - \uestII \vert))\right) - \int_{Z_{(j-1)}}^{Z_{(j)}} \exp(\psi^{(p)})(t) dt
\end{eqnarray*}
can be increased (and hence the log-likelihood) if $\psi^{(p)}_j$ were not diverging to $-\infty$. By a recursive reasoning, we obtain the same conclusion about the remaining integers.

\bigskip

\par \noindent Now suppose that $\lim_{p \to \infty} \max_{1\le i \le n} \psi^{(p)}(Z_{(i)})  = \infty$. Since $\psi^{(p)}$ is decreasing on $[0,\infty)$,  $\max_{1\le i \le n} \psi^{(p)}(Z_{(i)})  = \psi^{(p)}(Z_{(1)}) \to \infty $ as $p \to \infty$. For $x \in [Z_{(1)}, Z_{(2)}]$, we have
\begin{eqnarray*}
  1 &\ge & \int_{Z_{(1)}}^{Z_{(2)}} \exp(\psi^{(p)}_1) \exp\left( \frac{\psi^{(p)}_2 -\psi^{(p)}_1}{Z_{(2)} - Z_{(1)}} (t - Z_{(1)}) \right)  dt \\
  & = & \exp(\psi^{(p)}_1) (Z_{(2)} - Z_{(1)})  \frac{1- \exp( \psi^{(p)}_1 -\psi^{(p)}_2)}{\psi^{(p)}_1 -\psi^{(p)}_2} \\
  & \ge & \exp(\psi^{(p)}_1) (Z_{(2)} - Z_{(1)}) \frac{1}{\psi^{(p)}_1 -\psi^{(p)}_2 +1},
\end{eqnarray*}
where the last inequality follows from the fact that $\frac{1-\exp(-t)}{t} \ge \frac{1}{t+1}$ on $[0, \infty)$. This implies that
\begin{eqnarray*}
  \psi^{(p)}_1 -\psi^{(p)}_2   \ge \exp(\psi^{(p)}_1) (Z_{(2)} - Z_{(1)}) -1.
\end{eqnarray*}
Thus
\begin{eqnarray*}
  \psi^{(p)}_1 +\psi^{(p)}_2 = 2\psi^{(p)}_1 + \psi^{(p)}_2 -\psi^{(p)}_1 \le 2\psi^{(p)}_1 -\exp(\psi^{(p)}_1) (Z_{(2)} - Z_{(1)}) +1 \to -\infty
\end{eqnarray*}
as $p \to \infty$,  and hence $\psi^{(p)}_2 \le  (\psi^{(p)}_1 +\psi^{(p)}_2 )/2 \to -\infty$.  %
Since a maximizer of $\Phi_n$ is necessarily constant on $[-\tau_1, \tau_1]$ with $\tau_1$ the first kink point, here equal to $Z_{(1)}$, we can write
\begin{eqnarray*}
  1 = 2  \exp(\psi^{(p)}_1) Z_{(1)} + \sum_{j=2}^n \int_{Z_{(j)}}^{Z_{(j+1)}} \exp(\psi^{(p)}(t)) dt \ge 2  \exp(\psi^{(p)}_1) Z_{(1)}
\end{eqnarray*}
But this implies using the results obtained above that $\exp(\psi^{(p)}_1) \le 1/(2 Z_{(1)})$,  which contradicts the assumption that $\psi^{(p)}_1 \to \infty$.
\hfill $\Box$

\bigskip

\par \noindent \textbf{Proof of Proposition~\ref{prop:necessary-condition}.}
This follows from arguments similar to those used in the proof of Theorem 2.2
of \cite{duembgen_09}. \hfill $\Box$ %

\bigskip

\par \noindent \textbf{Proof of Proposition~\ref{AltChar}:} \ Using symmetry of $\widehat f_n$ and
$\Delta$, and the definition of $\widehat{\mathbb F}_n$, the inequality in
(\ref{CharMLE}) can be re-written as
\begin{eqnarray}\label{Char2}
  \int_{0}^\infty \Delta(x) d\widehat{\mathbb F}_n(x)  \le \int_{0}^\infty  \Delta(x) \widehat{f}^+_n(x)  dx,
\end{eqnarray}
 with equality for $\Delta$ satisfying $\widehat \psi^+_n \pm \epsilon
  \Delta \in \mathcal{SC}$ for $\epsilon > 0$ small enough.  For $z \in
[0, \infty)$, consider the concave perturbation function defined on $[0,
\infty)$ as
\begin{eqnarray*}
  \Delta_z(x) & = &  - (x - z)_+ \\
  & = &
  \begin{cases}
    z- x  \ \text{if $ x \ge z$} \\
    0 \ \ \text{otherwise}.
  \end{cases}
\end{eqnarray*}
Using Fubini's theorem, the inequality in (\ref{Char2}) yields
\begin{eqnarray*}
  Z_{(2n)} - z - \int_{z}^{Z_{(2n)}} \widehat{\mathbb F}_n(x) dx \ge Z_{(2n)} - z - \int_{z}^{Z_{(2n)}} \widehat F^+_n(x) dx.
\end{eqnarray*}
This yields
\begin{eqnarray*}
  \int_{0}^{z} \widehat F^+_n(x) dx \le \int_{0}^{z}   \widehat{\mathbb F}_n(x) dx + \int_{0}^{Z_{(2n)}} \widehat F^+_n(x) dx  - \int_{0}^{Z_{(2n)}}   \widehat{\mathbb F}_n(x) dx.
\end{eqnarray*}
Since for $z=0$, the inequality in (\ref{Char2}) becomes an equality
($\Delta_0$ is a straight line, and hence satisfies $\widehat \psi^+_n \pm
  \epsilon \Delta_0 \in \mathcal{SC}$ for $\epsilon > 0$ small enough), we
have that
\begin{eqnarray*}
  \int_{0}^{Z_{(2n)}} \widehat{\mathbb F}_n(x) dx = \int_{0}^{Z_{(2n)}} \widehat F^+_n(x) dx
\end{eqnarray*}
from which follows the inequality part in (\ref{CharMLE2}). The equality part
follows from noting that when $z$ is a knot of $ \widehat \psi^+_n =
\log(\widehat f^+_n)$ we have $\widehat \psi^+_n \pm \epsilon \Delta_z \in
\mathcal{SC}$ for $\epsilon > 0$ small enough.  \hfill $\Box$

\bigskip

\medskip

\bigskip

\bigskip

Next we state and prove propositions and lemmas needed for the proofs of  the
theorems in  Section~\ref{consistency}, and
then prove those theorems.
We start by stating several propositions of the  Glivenko-Cantelli type.  To
begin, let $\mathcal C$ be the class of convex and compact subsets
(intervals) in $\mathbb R$.  The following theorem of \cite{MR0436272} will
be an important component of the proof of Proposition~\ref{InfModul}
\begin{theorem}[Theorem 1.11, page 22, \cite{MR0436272}] %
  \label{GC}
  Let $P$ be a probability measure on $\mathbb R$.
  Then, as $n \to \infty$
  \begin{eqnarray*}
    \sup_{C \in \mathcal C} \vert \mathbb P_n(C)  - P(C) \vert \to 0
  \end{eqnarray*}
almost surely.
\end{theorem}

\medskip

\medskip

\begin{prop}\label{InfModul}
  Fix $\alpha \in (0,1)$. Let $X_1, \ldots, X_n$ be $n$ independent
  observations from the density $g^0 = \piest^0 f^0(\cdot - u^0_1) + (1- \piest^0)
  f^0(\cdot - u^0_2)$, and $\uestI$ and $\uestII$ be any consistent
  estimators of $u^0_1$ and $u^0_2$.  Write $\bar{u}_n = (\uestI +
  \uestII)/2, \bar{u}^0 = (u^0_1 + u^0_2)/2$ and let $\bar{N} = \sum_{i=1}^n
  1_{X_i \le \bar{u}}$. Then, for any $\epsilon > 0$, there exists $n_0 $
  such that
  \begin{eqnarray*}
    P\left( \min \bigg\{ \vert C \vert: C \ \in \mathcal C \ \textrm{such that} \ \bar{N}^{-1} \sum_{i: X_i \le \bar{u}_n} 1_{\{X_i \in C\}} \ge \alpha  \bigg \}  \ge  \frac{\alpha  G(\bar{u}^0)}{2 \Vert g^0 \Vert_\infty} \right) \ge 1- \epsilon
  \end{eqnarray*}
  for all $n \ge n_0$.
\end{prop}

\par \noindent
\textbf{Proof:} \
Let $\mathcal
C$ denote the class of all convex and compact subsets in $\mathbb R$. Fix
$\epsilon > 0$. For $C \in \mathcal C$ we have that
\begin{eqnarray*}
  \bar{N}^{-1} \sum_{i: X_i \le \bar{u}} 1_{\{X_i \in C\}} = \frac{\mathbb P_n(C \cap ]-\infty, \bar{u}])}{\mathbb G_n(\bar{u})} = \frac{\mathbb P_n(C_{\bar u})}{\mathbb G_n( \bar{u})}.
\end{eqnarray*}
where $C_{\bar u} = C \cap ]-\infty, \bar{u}])$ is also compact and convex.
To simplify notation, we write $C$ for $C_{\bar u_n})$, and $g$, $G$ and $P$
the true density, distribution function and associated measure respectively. Then
\begin{eqnarray*}
  \frac{\mathbb P_n(C)}{\mathbb G_n(\bar{u}_n)} - \frac{P(C)}{G(\bar{u}^0)} &= &  \frac{[\mathbb P_n(C) - P(C)] G(\bar u_n) - P(C) [\mathbb G_n(\bar u_n) - G(\bar u_n)]}{\mathbb G_n(\bar u_n) G(\bar u_n )}  \\
  && - \ \frac{P(C) [G(\bar u_n) - G(\bar u^0)]}{G(\bar u_n) G(\bar u^0)}.
\end{eqnarray*}
Consistency of $\uestI$ and $\uestII$, the fact that the class $\{]-\infty,
x]: x \in \RR\}$ is Glivenko-Cantelli (see Example 2.5.4 of
\cite{vdvwellner_96}) and continuity of $g$ imply that the event
$$\left \{G(\bar u_n) \ge \frac{G(\bar u^0)}{2}, \mathbb G_n(\bar u_n)  \ge
  \frac{G(\bar u^0)}{2},  \max_{t \in [\bar u^0 \wedge \bar u_n, \bar u^0 \vee
    \bar u_n]} g(t) \le 2g(\bar u^0) \right\}$$
 occurs with probability converging to one. On the other hand, we have that
\begin{eqnarray*}
  \left \vert \frac{\mathbb P_n(C)}{\mathbb G_n(\bar{u}_n)} - \frac{P(C)}{G(\bar{u}^0)}  \right \vert &\le &   \frac{\sup_{C \in \mathcal C} \vert \mathbb P_n(C) - P(C) \vert + \sup_{x \in \RR} \vert \mathbb G_n(x) - G(x) \vert}{\mathbb G_n(\bar u_n) G(\bar u_n )}   \\
  && \ + \frac{\vert G(\bar u_n) - G(\bar u^0)\vert }{G(\bar u_n) G(\bar u^0)} \\
  & \le & \frac{\sup_{C \in \mathcal C} \vert \mathbb P_n(C) - P(C) \vert + \sup_{x \in \RR} \vert \mathbb G_n(x) - G(x) \vert}{\mathbb G_n(\bar u_n) G(\bar u_n )}   \\
  && \ + \frac{g(\xi) \vert \bar u_n - \bar u^0 \vert }{G(\bar u_n) G(\bar u^0)}, \ \ \textrm{for some $\xi \in [ \bar u^0 \wedge \bar u_n, \bar u^0 \vee \bar u_n ]$} \\
  & \le & D \left( \sup_{C \in \mathcal C} \vert \mathbb P_n(C) - P(C) \vert + \sup_{x \in \RR} \vert \mathbb G_n(x) - G(x) \vert + \vert \bar u_n - \bar u^0 \vert \right)
\end{eqnarray*}
probability increasing to one, where $ D = 4 G(\bar u^0)^{-2} \max(2 g(\bar
u^0),1)$. Using
Theorem~\ref{GC},
and again consistency of $\bar u_n$, we conclude that
\begin{eqnarray*}
  \sup_{C \in \mathcal C} \left \vert \frac{\mathbb P_n(C)}{\mathbb
      G_n(\bar{u}_n)} - \frac{P(C)}{G(\bar{u}^0)}  \right \vert \to 0  \  \
  \textrm{in probability as } n \to \infty.
\end{eqnarray*}
In particular, this implies that for any fixed $\epsilon > 0$ there exists
$n_0$
such that
\begin{eqnarray*}
  P\left(\frac{\mathbb P_n(C)}{\mathbb G_n(\bar{u}_n)} \le \frac{P(C)}{G(\bar{u}^0)} + \frac{\alpha}{2} \ \textrm{for all $C \in \mathcal C$} \right) \ge 1-\epsilon.
\end{eqnarray*}
Let us denote by $\mathcal E_{n,\alpha/2}$ the event   $\left \{\frac{\mathbb P_n(C)}{\mathbb G_n(\bar{u}_n)} \le \frac{P(C)}{G(\bar{u}^0)} + \frac{\alpha}{2} \ \textrm{for all $C \in \mathcal C$} \right \}$. Also, consider the subclass
\begin{eqnarray*}
  \mathcal C_{\alpha/2} = \left\{C \in \mathcal C: \frac{P(C)}{G(\bar{u}^0)} \ge \frac{\alpha}{2}  \right \}.  %
\end{eqnarray*}
Clearly, for any $C \in \mathcal C_{\alpha/2}$ we have that
$$\vert C \vert \ge \frac{\alpha G(\bar{u}^0)}{2 \Vert g \Vert_\infty}.$$
On the other hand,
\begin{eqnarray*}
  && \left \{\ C \in \mathcal C: \bar{N}^{-1} \sum_{i: X_i \le \bar{u}_n} 1_{\{X_i \in C\}} \ge \alpha \right\}  =  \left\{\ C \in \mathcal C: \frac{\mathbb P_n(C)}{\mathbb G_n(\bar{u}_n)} \ge \alpha \ \textrm{for $C \in \mathcal C$}  \right \} \\
  && =  \left\{\left\{C \in \mathcal C: \frac{\mathbb P_n(C)}{\mathbb G_n(\bar{u}_n)} \ge \alpha \right \} \cap \mathcal{E}_{n,\alpha/2}\right\} \  \cup \   \left\{\left\{C \in \mathcal C: \frac{\mathbb P_n(C)}{\mathbb G_n(\bar{u}_n)} \ge \alpha \right \} \cap \mathcal{E}^c_{n,\alpha/2}\right \} \\
  && \subseteq  \mathcal C_{\alpha/2} \   \cup \  \left\{\left\{C \in \mathcal C: \frac{\mathbb P_n(C)}{\mathbb G_n(\bar{u}_n)} \ge \alpha\} \right \} \cap \mathcal{E}^c_{n,\alpha/2}\right \}
\end{eqnarray*}
implying that for $n \ge n_0$
\begin{eqnarray*}
  P\left(\min \big\{ \vert C \vert: C \ \in \mathcal C \ \textrm{such that} \ \bar{N}^{-1} \sum_{i: X_i \le \bar{u}_n} 1_{\{X_i \in C\}} \ge \alpha  \big \}  <  \frac{\alpha G(\bar u^0)}{2 \Vert g^0 \Vert_\infty} \right)  \le \epsilon
\end{eqnarray*}
completing the proof of the Proposition. \hfill $\Box$

\medskip

\bigskip

\begin{lemma}\label{KL1}
  Let $\psi \in \mathcal{SC}$,  $a < b \ \in \RR$ and $\lambda \in [0,1]$. If $p(t) = \lambda \exp(\psi(t - a))  + (1- \lambda) \exp(\psi(t - b))$, then
  \begin{eqnarray*}
    \int_{\RR} \vert \log(p(t)) \vert  p(t)  dt  < \infty.
  \end{eqnarray*}
\end{lemma}

\par \noindent  %
\textbf{Proof:}
Note that for $t \in \RR$
we have that
\begin{eqnarray*}\label{Res}
  \vert t - a \vert \le \vert t - b \vert \Longleftrightarrow  t \le (a+b)/2.
\end{eqnarray*}
Using the fact that the function $ t \mapsto \vert \log (t) \vert t $ is
increasing on $(0, e^{-1}]$ and that $\psi$ is decreasing on $[0, \infty)$,
we have for $A$ taken to be larger than $(a+b)/2$,
\begin{eqnarray*}
  \int_{A}^ \infty \Big \vert \log(p(t)) \vert p(t) dt  &\le & \int_A^\infty \vert \psi(t - a) \vert  \exp(\psi(t-a)) dt  \\
  & < & \infty
\end{eqnarray*}
 where the second inequality follows e.g.\ from Lemma 2 in \cite{schuhmacher_10}.  By symmetry, the integral on $(-\infty, -A]$ is also finite.   \hfill $\Box$

\bigskip

\par \noindent In the following, we introduce a condition that will
be useful in proving consistency of our estimator.

\bigskip

\noindent \textbf{The condition (C):}\ We will say that an (ordered) pair
$(f,h) = (e^\psi, e^\phi)$ of log-concave densities satisfy the condition
$(C)$ if for all $d \in \RR$,
\begin{eqnarray*}
  \int_0^\infty \left \vert \phi(t +d) \right \vert f(t) dt
  < \infty.  %
\end{eqnarray*}
Note that $h$ must have support $(-\infty,\infty)$.

\medskip

\begin{lemma}\label{KL2}

  Let $\psi_0$ and $\phi$ be elements of $\mathcal{SC}$.  Let $a_0 < b_0$ be
  real numbers and $\lambda_0 \in [0,1]$.  Let $(a_k)$, $(b_k)$, and
  $(\lambda_k)$, for $k > 0$, be sequences of real numbers converging to
  $a_0, b_0 $, and $\lambda_0$, respectively.  Let
  \begin{align*}
    p_k( \cdot ) & = \lambda_k \exp(\psi_0(  \cdot -a_k)) + (1-\lambda_k)
    \exp(\psi_0( \cdot -b_k)),\\
    p_{k,\phi}( \cdot ) & = p_k( \cdot ) \star e^{\phi(\cdot)}
    = \lambda_k \exp(\psi_0(\cdot -a_k)) \star e^{\phi(\cdot)} + (1-\lambda_k)
    \exp(\psi_0(t-b_k)) \star e^{\phi(\cdot)},
  \end{align*}
  for $k=0, 1, \ldots$.  Let $e^{\psi_0} $ and $e^\phi$ satisfy condition
  $(C)$.  Then
  \begin{equation}
    \int_{\RR} \log \lp p_{k,\phi}(t) \rp p_0(t) \, dt \longrightarrow
    \int_{\RR} \log \lp p_{0,\phi}(t) \rp p_0(t) \, dt
    \label{eq:integral-convergence}
  \end{equation}
  as $k \to \infty$.
\end{lemma}

\bigskip

\par \noindent
\textbf{Proof:} \
It is enough to show that for $k$ large, $|
\log p_{k,\phi} | $ is bounded by a $p_0$-integrable function to apply the
Lebesgue Dominated Convergence Theorem.  Since $p_{k,\phi}$ is uniformly
bounded above, to upper bound $| \log p_{k, \phi}|$ we need to lower bound
$p_{k,\phi}$. Now,
  \begin{align}
    p_{k,\phi}(x)
    & = \int_{\RR} p_{k}(y) e^{\phi(x-y)} \, dy
    \ge \int_{\alpha}^\beta  \lp  \min_{z \in [\alpha,\beta]}    p_{k}(z) \rp
    e^{\phi(x-y)}      \, dy
    \label{eq:bound1}
  \end{align}
  for any $\alpha$ and $\beta$.  Taking $\alpha $ and $\beta$ to be in the
  interior of the support of $p_0$, we have that, for all $k$ large enough,
  $\min_{z \in [\alpha,\beta]} p_{k}(z) \ge c $ for some $c>0$.  Now because
  $\phi$ is concave, if $x \ge A$ where $A$ is large enough, $\phi$ is
  decreasing on $[x-\beta, x-\alpha]$, i.e.\ ${\phi(x-y)}$ is increasing in
  $y$ for $y \in [\alpha,\beta]$, so that then
  $\phi(x-y) \ge \phi(x-\alpha)$.  Taking $A$ large enough that
  $\phi(x-\alpha) < 0 $, so $|\phi(x-y) | \le |\phi(x-\alpha)|$, we have that,
  for $k$ large, and $x > A$, \eqref{eq:bound1} is bounded below by
  \begin{equation*}
    c \int_\alpha^\beta e^{\phi(x-\alpha)} \, dy
    = (\beta - \alpha) c e^{\phi(x-\alpha)}.
  \end{equation*}
  Recall that if $d_1 < d_2$, then
  \begin{align*}
    |t-d_1| \le | t-d_2|
    \Longleftrightarrow
    t \le (d_1+d_2)/2.
  \end{align*}
  Thus, taking $A > (a+b)/2$,
  \begin{align*}
    \int_A^\infty |\phi(x-\alpha)| p_0(x) \, dx
    \le \int_A^\infty |\phi(x-\alpha) | e^{\psi_0(x-b_0)} \, dx
    = \int_{A+b}^\infty |\phi(x-\alpha+b)| e^{\psi_0(x)}\, dx,
  \end{align*}
  and this is less than infinity, by the condition $(C)$.
  Thus, for $x \ge A$ and $k$ large, %
  \begin{align*}
    |\log p_{k,\phi}(x) |
    \le |\log \lp (\beta-\alpha)c\rp | + |\phi(x-\alpha)|,
  \end{align*}
  so we have shown that $|\log p_{k,\phi}|$ is bounded by a function
  $p_0$-integrable on $[A,\infty)$.  A similar argument shows it is bounded
  by a function integrable on $(-\infty,-A]$.  Finally, since $p_{0,\phi}$
  has support $(-\infty,\infty)$ (since $e^\phi$ has support
  $(-\infty,\infty)$ by the condition $(C)$),
  \begin{align*}
    \sup_{x \in [-A,A]} |\log p_{k,\phi}(x) | \le \sup_{x \in [-A-\eta,
      A+\eta]} |\log p_0(x)|
    < \infty,
  \end{align*}
  for all  $k$ large enough that $\max(|a_0-a_k|, |b_0-b_k|) < \eta$ for
  some $\eta > 0$.  Thus for $k$ large we have shown that $|\log p_{k,\phi}|
  $ is bounded by a $p_0$-integrable function on $\RR$, and $p_{k,\phi}(t) \to
  p_{0,\phi}(t) $ for all $ t \in \RR$, so \eqref{eq:integral-convergence}
  holds.  \hfill $\Box$

\bigskip

Note that for any log-concave density $f = e^\psi$, if $\phi$ is a polynomial
function, then $f$ and $e^\phi$ satisfy condition (C), since log-concave
densities have finite moments of all orders.  In particular, we may take
$e^\phi$ to be any normal density, which we will do below in the proof of
consistency.  Theorem~\ref{ConvMixed} will rely on the following two
Glivenko-Cantelli propositions.  As is standard, we say that a class of
functions $\cal F$ is a {\em $G$-Glivenko-Cantelli class} if $\sup_{f \in
  \mathcal{F}} | \GG_n(f) - G(f) | \to 0$ as $n \to \infty$ almost surely,
where
$\GG_n$ is the
empirical probability measure based on i.i.d.\ draws from $G$, and we write
$\mu(f) = \int f d\mu$ for a function $f$ and a measure $\mu$
(as in e.g.\ \cite{vdvwellner_96} or \cite{vandeGeer:2009tk}).
(This is sometimes known as being {\em strong} Glivenko-Cantelli because the
convergence is almost sure, rather than in probability \citep{vanderVaart:2000vr}.)

\bigskip

\begin{prop}\label{Glivenko1}
  Let $M> 0$ and $\mathcal{LC}_M$ be the class of log-concave functions
  bounded by $M$. For a fixed $b > 0$, consider the class of functions
  \begin{eqnarray*}
    \mathcal F_b = \Big \{ f:  f= \log(\lambda h_1 + (1-\lambda) h_2 + b), \ h_1, h_2  \in \mathcal{LC}_M, \ \lambda \in [0,1] \Big \}.
  \end{eqnarray*}
  Then $\mathcal F_b$ is a Glivenko-Cantelli class for any probability measure $P$.
\end{prop}

\par \noindent %
\textbf{Proof:} \
Since all $f
\in {\cal LC}_M$ are unimodal, ${\cal LC}_M \subset {\cal M}_M^+
- {\cal M}_M^+$ where ${\cal M}_M^+ := \lb f: 0 \le f \le M, f \mbox{ is
  nondecreasing} \rb$.  By Theorem 2.7.5 and Theorem 2.4.1 in
\cite{vdvwellner_96}, ${\cal M}_M^+$ is a (strong) $P$-Glivenko-Cantelli class
for any $P$.
By the Glivenko-Cantelli preservation Theorem 3 of \cite{vanderVaart:2000vr},
$\varphi({\cal M}_M^+, {\cal M}_M^+)$ is $P$-Glivenko-Cantelli, with
$\varphi(x,y) := x-y$ here, so ${\cal LC}_M$ also is $P$-Glivenko-Cantelli, since
it has an integrable envelope by construction.
Now letting
 $\varphi(h_1, h_2, \lambda) := \log \lp \lambda h_1 + (1-\lambda ) h_2 +
b \rp$, another application of
 Theorem 3 of \cite{vanderVaart:2000vr} yields that
 ${\cal F}_b = \varphi \lp {\cal LC}_M, {\cal LC}_M,
 \lb \lambda  \mathbbm{1}_{(-\infty,\infty)}(x) : \lambda \in [0,1] \rb \rp$
is also Glivenko-Cantelli since each of the three argument classes are,
as long as ${\cal F}_b $ has a $P$-integrable envelope.  But this is immediate
since $\log(b) \le f \le \log(M+b)$ for all $ f \in {\cal F}_b$. \hfill $\Box$

\bigskip

We will also need a Glivenko-Cantelli result for the case $b=0$.  In order
for such a result to hold, we need to restrict attention to just logs of
mixtures of a single density; despite the unboundedness of the functions
involved, this will allow us to have an integrable envelope, at least for the
true $P_0$.

\begin{prop}\label{Glivenko2}
  Fix $\tilde f_0 \in \mathcal{SLC}_1$, $\piest \in [0,1]$, $a < b $ and let
  $p_0(x) = \piest \tilde f_0(x-a) + (1-\piest) \tilde f_0(x-b), \ x \in \RR $, and
  $P_0$ the corresponding probability measure. Let $q \in \mathcal{SLC}_1$ be
  such that $(\tilde f_0, q)$ satisfy the Condition (C).
  Then the class of functions
  \begin{align*}
    \mathcal H &
    =  \Big \{ \log \left(\lambda q (x- \alpha) +
      (1- \lambda) q (x-\beta) \right) :  \\
    &  \hspace{1cm} \ \alpha \in \big[a- r, a +  r \big], \beta
    \in \big[b - r, b+ r \big], \lambda \in [0,1] \Big \}
  \end{align*}
  for any fixed $r >0$
  is $P_0$-Glivenko-Cantelli.
\end{prop}

\par \noindent %
\textbf{Proof:} \
 Without loss of
generality, we take $r=1$ in this proof. For any $c \in \RR$ let
$\mathcal{P}_{c} := \lb q (\cdot - z) : z \in [c-1,c+1] \rb$.  Then the
classes ${\cal P}_{a}$ and ${\cal P}_{b}$ are both contained in
$\mathcal{LC}_M$ with $M = \sup_x q(x)$, and in the proof of
Proposition~\ref{Glivenko1} we showed that $\mathcal{LC}_M$ is
$P$-Glivenko-Cantelli for any $P$. Let $T(h_1, h_2, \lambda) := \log \lp
\lambda h_1 + (1-\lambda) h_2 \rp$ and let $ I_{0,1}: = \lb \lambda
\mathbbm{1}_{(-\infty,\infty)}(x) : \lambda \in [0,1] \rb$.  Then by Theorem
3 of \cite{vanderVaart:2000vr}, the class $T \lp \mathcal{P}_{a},
\mathcal{P}_{b}, I_{0,1} \rp$ is $P_0$-Glivenko-Cantelli as long as it has a
$P_0$-integrable envelope, which is what we will show now.

Let $c_+ \in \RR$ be such that $\log q(x - (b + 1)) < 0$ if $x > c_+$. Then
for any $h \in \mathcal{H}$, $h(x) \le 0 $ for $x > c_+$ by the unimodality
of $q$.  Note $q(x - \alpha ) \ge q(x - (a-1))$ for $x \ge a+ 1 \ge \alpha
\ge a-1$ and $q(x-\beta) \ge q(x - (b-1))$ if $x \ge b+1 \ge \beta \ge b-1 $.
Thus, for $x > (b+1) \vee c_+$, $h \in \mathcal{H}$ satisfies
\begin{align}
  0 \ge h(x)
  & \ge \log \lp \lambda q(x - (a-1)) + (1-\lambda) q
  (x-(b-1)) \rp \nonumber \\
  & \ge \log q(x - (a-1)), \label{temp}
\end{align}
since $b > a$.  A similar argument shows for $x < (a -1) \wedge c_-$ for some
$c_- \in \RR$, that $0 \ge h(x) \ge \log q(x - ( b + 1))$.  By the assumption
that Condition $(C)$ holds for $(\tilde f_0, q)$, it is immediate that $\log
q(\cdot - (a-1))$ and $\log q(\cdot - (b+1))$ are $p_0$-integrable.  All $h
\in \mathcal{H}$ are uniformly bounded on $[(a-1) \wedge c_-, (b+1) \vee
c_+]$, since the support of $p_{0,\varphi}$ is $\RR$, since $p_{0,\varphi}$
is unimodal, and since the range of $\alpha$ and $\beta$ in the definition of
$\cal H$ is bounded.  Thus $\mathcal{H}$ has a $P_0$-integrable envelope and
so we are done. \hfill $\Box$

\bigskip

\medskip

\begin{theorem}\label{Boundedness}
  Let $\widehat g_n$ denote again the MLE of the mixed density. Then, for a
  given $\epsilon > 0$ there exist an integer $n_0 > 0$ and $D > 0$ depending
  on $\epsilon$ such that
  \begin{eqnarray*}
    P\left(\Vert \widehat g_n \Vert_\infty \le D \right) \ge  1 - \epsilon.
  \end{eqnarray*}
  for all $n \ge n_0$.
\end{theorem}

\bigskip

\bigskip

\par \noindent
\textbf{Proof:} \
Fix $\epsilon > 0$.  Recall  that if $a < b$, then for $t \in \RR$
\begin{eqnarray*}\label{Res}
  \vert t - a \vert \le \vert t - b \vert \Longleftrightarrow  t \le (a+b)/2.
\end{eqnarray*}
Write $\bar{u} = (\uestI + \uestII)/2$. Using the symmetry of $\widehat{\psi}_n$ and the fact that it is decreasing on $[0, \infty)$, it follows that
\begin{eqnarray}\label{BoundLik}
  \Phi_n(\log(\widehat g_n))  \le \frac{1}{n} \sum_{i: X_i \le \bar{u}} \widehat \psi_n(X_i - \uestI) + \frac{1}{n} \sum_{i: X_i > \bar{u} } \widehat \psi_n(X_i - \uestII) -1.
\end{eqnarray}
Without loss of generality, we can assume that $\{X_{1}, X_{2},  \ldots, X_{\bar{N}}\} = \{X_i: X_i \le \bar{u} \}$ and that  $\psi_n(X_{1} - \uestI) < \ldots < \widehat \psi_n(X_{\bar N} - \uestI)$,and that
\begin{eqnarray*}
  \max_{i: X_i > \bar u} \widehat \psi_n(X_i - \uestII) \le \max_{i: X_i \le \bar u} \widehat \psi_n(X_i - \uestI)
\end{eqnarray*}
otherwise the same reasoning above has to be applied to the second term of the right side of  the inequality in (\ref{BoundLik}).  In the following, let
$$
M_n = \widehat \psi_n(X_n - \uestI)
$$
the largest value taken by $\widehat \psi_n(X_i - \uestI)$ and $\widehat \psi_n(X_i - \uestII)$ for $i =1, \ldots, n$.

\medskip

\par \noindent Following the notation of \cite{schuhmacher_10} in the proof of their Lemma 4, let $k = \lfloor \bar{N} \alpha \rfloor $ for some fixed $\alpha \in (0,1]$, $\mathcal C_{k}$ the convex hull of the set $\{X_i -\uestI, i \ge k \}$ and $m = \widehat \psi_n( X_k - \uestI )$.  Also, let
$$\eta = \frac{\alpha G(\bar u^0)}{2 \Vert g^0 \Vert_\infty}.$$

\smallskip

\par \noindent In the sequel,  $\mathcal C_k + \uestI = \{t +\uestI: t \in \mathcal C_k \}$.  Note that $\mathcal C_k + \uestI$ is convex and compact such that
\begin{eqnarray*}
  \sum_{i: X_i \le \bar{u}} 1_{\{X_i \in \mathcal C_k + \uestI\}} = \sum_{i=1}^{\bar N} 1_{\{X_i \in \mathcal C_k + \uestI\}} = \bar N - k +1 \ge \bar N (1-\alpha).
\end{eqnarray*}
By taking $1-\alpha$ in place of $\alpha$ in Proposition \ref{InfModul}, we can find an integer $n_0 > 0$ such that for $n \ge n_0$.
\begin{eqnarray}\label{Bound}
  P\left(\vert C_k \vert \ge \eta \right) \ge  1-\epsilon/5.
\end{eqnarray}
Consider the event
\begin{eqnarray}\label{IneqMn}
  \Big\{M_n > B_0: = \max(\log(1/\eta) + 1 , 0) \Big \}.
\end{eqnarray}
Using the same argument as in \cite{schuhmacher_10}, we can write
\begin{eqnarray}\label{Ineqm}
  1= \int_{\RR} e^{\widehat \psi_n(t)}  dt &\ge & \int_{\mathcal C_{k}} e^{\widehat \psi_n(t)} dt  \nonumber\\
  & \ge & e^m \vert \mathcal C_{k} \vert
\end{eqnarray}
using monotonicity of $\widehat \psi_n$ on the convex set $\mathcal C_k$. By (\ref{IneqMn}). Hence,
\begin{eqnarray}\label{Event0}
  \Big \{ M_n > B_0, \vert \mathcal C_k \vert > \eta \Big \} \subset \{ M_n > m+1 \}.
\end{eqnarray}
Indeed, this follows from the fact that
\begin{eqnarray*}
  e^{M_n - 1} & > & 1/\eta, \ \ \ \textrm{by (\ref{IneqMn})} \\
  & > & 1/\vert \mathcal C_k \vert  \\
  & \ge & e^m \ \textrm{by (\ref{Ineqm})}.
\end{eqnarray*}
Also, on the event on the right side of (\ref{Event0}), we have for $t$ such that $t \in \mathcal C_k + \uestI$
\begin{eqnarray*}
  \widehat \psi_n\left( X_n-\uestI + \frac{t- X_n}{M_n-m} \right)   & \ge & \frac{1}{M_n-m} \widehat \psi_n(t -\uestI)+ \frac{M_n-m-1}{M_n-m}\widehat \psi_n(X_n-\uestI)  \\
  & \ge & \frac{m}{M_n-m} + \frac{M_n(M_n-m-1)}{M_n-m}  = M_n-1.
\end{eqnarray*}
Thus,
\begin{eqnarray*}
  \textrm{Leb}(x \in \RR: \widehat \psi_n(x) \ge M_n-1) \ge \textrm{Leb}\left( \frac{M_n-m-1}{M_n-m} (X_n-\uestI) + \frac{1}{M_n-m} \mathcal{C}_k   \right)  = \frac{\vert \mathcal C_k \vert }{M_n-m},
\end{eqnarray*}
and hence
$$ 1 = \int_{\RR} e^{\widehat \psi_n(t)}  dt \ge \frac{e^{M_n-1} \vert \mathcal C_k \vert}{M_n- m} $$
Hence,  the event on the right side of (\ref{Event0}) is included in
\begin{eqnarray*}
  \Big \{ m \le M_n - e^{M_n -1} \vert \mathcal C_k \vert < M_n - e^{M_n -1} \eta \Big \}.
\end{eqnarray*}
Let $B_1 > 0$ be large enough so that $B_1 > B_0$,  $B_1 - e^{B_1 - 1} \eta < -e^{B_1 - 1} \eta/2$ and $x \mapsto x - e^{x-1} \eta/2 $ is decreasing on $[B_1, \infty)$, and consider the event
\begin{eqnarray*}\label{Event1}
  \Big \{ M_n > B_1, \vert \mathcal{C}_k \vert > \eta \Big \}.
\end{eqnarray*}
Then, this event is included in
\begin{eqnarray*}
  \Big\{  m \le -e^{M_n - 1} \eta/2 \Big \}
\end{eqnarray*}
and occurrence of the latter implies in turn that
\begin{eqnarray*}
  \Phi_n(\log(\widehat g_n)) \le \frac{1}{n} \left[(\lfloor \alpha \bar N  \rfloor- 1) M_n - (\bar N - \lfloor \alpha \bar N \rfloor +1 ) \frac{1}{2} e^{M_n-1} \eta + \frac{n - \bar N}{n} M_n \right ],
\end{eqnarray*}
and consequently
\begin{eqnarray*}
  \Phi_n(\log(\widehat g_n)) \le \left(1 - (1-\alpha) \frac{\bar N}{n}\right) M_n - \frac{(1-\alpha) \bar{N} }{2n} e^{M_n-1} \eta.
\end{eqnarray*}
Using similar arguments as in the proof of our Proposition \ref{InfModul}, we can increase $n_0$ to ensure that for $n \ge n_0$
\begin{eqnarray*}
  P\left( \frac{\bar N}{n}  \ge \frac{G^0(\bar u^0)}{2} \right) \ge 1 - \epsilon/5.
\end{eqnarray*}
Hence, occurrence of the event
\begin{eqnarray*}\label{Event2}
  \Bigg \{ M_n > B_1, \vert C_k \vert > \eta, \frac{\bar N}{n}  \ge \frac{G^0(\bar u^0)}{2} \Bigg \}.
\end{eqnarray*}
implies that
\begin{eqnarray}\label{LowerBound}
  \Phi_n(\log(\widehat g_n)) \le \left(1 - (1-\alpha) \frac{G^0(\bar u^0)}{2}\right) M_n - \frac{(1-\alpha) G^0(\bar u^0) }{4} \eta e^{M_n-1}
\end{eqnarray}
Next, if we take $e^w = \varphi_\sigma$ the density of standard Normal random
variable.
Then, %
\begin{eqnarray*}
  \Phi_n(\log(\widehat g_n))  & \ge &  \Phi_n\left(\piest e^{w(\cdot - \uestI)} + (1-\piest) e^{w(\cdot - \uestII)}\right) \\
  & = & \int_{\RR} \log\left(\piest e^{w(x - \uestI)} + (1-\piest) e^{w(x - \uestII)}\right) d(\mathbb G_n(x) - G^0(x))\\
  && \ +  \int_{\RR} \log\left(\piest e^{w(x - \uestI)} + (1-\piest) e^{w(x- \uestII)}\right)  g^0(x) dx -1.
\end{eqnarray*}
For $\eta > 0$ small enough, it follows from consistency of $\uestI, \uestII$ and $\piest$ that the event
\begin{eqnarray*}
  \mathcal E^0 = \Big \{\vert \piest - \piest^0 \vert \le \eta/2, \vert \uestI - u^0_1 \vert \le \eta/2, \vert \uestII - u^0_2 \vert \le \eta/2 \Big\}
\end{eqnarray*}
occurs with probability greater than $1-\epsilon/5$ for $n \ge n_0$ (at the
cost of increasing the previous $n_0$).  By Proposition \ref{Glivenko2},
since $(f^0, e^w)$ certainly satisfy Condition $(C)$,
this implies that the event
\begin{eqnarray*}
  \mathcal E_{\text{emp}} = \Bigg \{\left \vert \int_{\RR} \log\left(\piest e^{w(x - \uestI)} + (1-\piest) e^{w(x - \uestII)}\right) d(\mathbb G_n(x) - G^0(x)) \right \vert \le \frac{1}{2} \Bigg \}
\end{eqnarray*}
occurs with probability greater than $1- \epsilon/5$ for $n \ge n_0$.  Also, there exists a constant $D > 0$ (depending on $\epsilon$)
such that the event
\begin{eqnarray*}
  \mathcal E_{KL} = \Bigg \{\left \vert \int_{\RR} \log\left(\piest e^{w(x - \uestI)} + (1-\piest) e^{w(x- \uestII)}\right) g^0(x) dx \right \vert \le D  \Bigg\}
\end{eqnarray*}
with probability greater than $1-\epsilon/5$.  Indeed, note that $\piest e^{w(x - \uestI)} + (1-\piest) e^{w(x- \uestII)} \le 1$ for all $x \in \RR$ so that
\begin{eqnarray*}
&&\left \vert \int_{\RR} \log\left(\piest e^{w(x - \uestI)} + (1-\piest) e^{w(x- \uestII)}\right) g^0(x) dx \right \vert \\
&&= - \int_{\RR} \log\left(\piest e^{w(x - \uestI)} + (1-\piest) e^{w(x- \uestII)}\right) g^0(x) dx.
\end{eqnarray*}
Using this fact and concavity of the logarithm, it follows that
\begin{eqnarray*}
0 &\le &  -\int_{\RR} \log\left(\piest e^{w(x - \uestI)} + (1-\piest) e^{w(x- \uestII)}\right) g^0(x) dx \\
& \le &  \frac{1}{2}  \int_{\RR} \left(\piest (x-\uestI)^2 + (1-\piest) (x-\uestII)^2\right) g^0(x)dx \\
& < & \infty
\end{eqnarray*}
since log-concave densities have finite moments of any order.   Put $L = -3/2 -D$, and let $B_2  > B_1$ such that
\begin{eqnarray*}
  \left(1 - (1-\alpha) \frac{G(\bar u^0)}{2}\right) B_2 - \frac{(1-\alpha) G(\bar u^0) }{4} \eta e^{B_2-1} < L
\end{eqnarray*}
and the function $x \mapsto \left(1 - (1-\alpha) \frac{G(\bar u^0)}{2}\right) x - \frac{(1-\alpha) G(\bar u^0) }{4} \eta e^{x-1}$ is decreasing on $[B_2, \infty)$. By the choice of $B_2$ above, we see that \begin{eqnarray*}
  \mathcal E_{\cap} = \Bigg\{ M_n > B_2, \vert C_k \vert > \eta, \frac{\bar N}{n}  \ge \frac{G^0(\bar u^0)}{2} \  \Bigg \} \cap \mathcal E^0 \cap \mathcal E_{\text{emp}} \cap \mathcal E_{KL} = \emptyset
\end{eqnarray*}
Hence,
\begin{eqnarray*}
  P\Big(\Vert \widehat g_n \Vert_\infty > B_2\Big) = P(M_n > B_2)  = P(M_n > B_2, \mathcal E_{\cap}^c) \le P(\mathcal E_{\cap}^c) \le \epsilon.
\end{eqnarray*}
Note that $B_2$ depends on $\epsilon$ through $L$.    \hfill $\Box$

\bigskip
\bigskip

To establish the first consistency result, note by the definition of
$\widehat g_n$ we have that for any density $\tilde g$ of the form
\begin{eqnarray*}
  \label{g}
  \tilde g = \piest f(\cdot - \uestI) + (1- \piest) f(\cdot - \uestII).
\end{eqnarray*}
for $f \in \mathcal{SLC}_1$ that
\begin{eqnarray*}
  0 \le \int \log \widehat g_n \, d\GGn - \int \log \tilde g \, d\GGn
  \le \int \log (\widehat g_n + b) \, d\GGn - \int \log \tilde g \, d\GGn.
\end{eqnarray*}
Now let $b > 0$, and $G^0$ be the cdf of the true mixed density $g^0$.  As
first established by \cite{pal_07} in their Lemma 1, the inequality above
yields
\begin{equation}
  \label{eq:Hellinger}
  \begin{split}
    2 H^2(\widehat g_n, g^0)
    & \le   \epsilon(b) + \int_{\RR} \log(\widehat g_n(t) + b)d(\mathbb G_n(t) - G^0(t))  \\
    & \quad + \ \int_{\RR} \log(g^0(t) +  b)dG^0(t) - \int_{\RR} \log(\tilde{g}(t)) dG^0(t)  \\
    & \quad  - \int_{\RR} \log(\tilde{g}(t)) d(\mathbb G_n(t) - G^0(t))
  \end{split}
\end{equation}
where $H(p,q)$ denotes the Hellinger distance and
\begin{eqnarray*}
  \epsilon(b) = 2 \int_{\RR} \sqrt{\frac{b}{b+ g^0(t)}}dG^0(t). %
\end{eqnarray*}
The main idea behind introducing the small positive quantity $b$ is to avoid
integration issues due to the fact that $\log(\widehat g_n)= -\infty$ outside
a certain interval. In the problem of estimating a log-concave density, the
empirical term on the right side of the previous inequality was shown to
converge to zero using the fact that the maximum value of the log-concave MLE
stays bounded in $n$, see Theorem 3.2 of \cite{pal_07} and Lemma 4 in
\cite{schuhmacher_10} for a generalization of the same result in the
multivariate setting. This property of the MLE was then combined with the
fact that a level set of a bounded unimodal function is convex and
compact. This cannot be claimed anymore if we replace a unimodal function by
a mixture of two unimodal functions, not even when those functions are
log-concave.  For this reason, we shall use instead the Glivenko-Cantelli
results proved above.

In our proof, the natural choice for $\tilde g$ would be
\begin{equation}
  \label{eq:defn:widecheck-g}
  \widecheck{g}_{n} = \piest f^0(\cdot - \uestI) + (1- \piest) f^0(\cdot - \uestII),
\end{equation}
but this is problematic.  For instance $\int \log \widecheck{g}_n d\GGn$ can be $-\infty$
for arbitrarily large $n$ if $g^0$ has compact support.  This might occur
when the smallest/ largest order statistic is very close to the left/right
end of the support if $f^0$ and $\uestI$/$\uestII$ is smaller/larger than
$\uestI^0$/$u^0_2$.  But even if $g^0$ does not have compact support, if it does
have too large a slope (e.g., it approximates having compact support by
dropping towards $0$ very quickly at some point), then  $\int \log \widecheck{g}_n \,
dG^0$ can still be infinite.  Thus, we consider the surrogate density function
\begin{eqnarray}\label{gh}
  \widecheck{g}_{n,h} = \piest f^0_h(\cdot - \uestI) + (1- \piest) f^0_h(\cdot - \uestII)
\end{eqnarray}
where $f^0_h = f^0 \star \phi_h$, the convolution of $f^0$ and the density
of centered normal with standard deviation $h > 0$. This choice alleviates
the above-mentioned problems; the risk of having a divergent log-likelihood
is now excluded since $f_h^0$, the component density of
$\widecheck{g}_{n,h}$,
is supported on
$\RR$.  Note that $f^0_h$ is a log-concave density by
preservation of log-concavity under convolution \citep{Ibragimov56}, and is
symmetric.

\bigskip
\par \noindent \textbf{Proof of Theorem~\ref{ConvMixed}:} \ The proof starts
from the inequality in \eqref{eq:Hellinger} with $\widecheck g_{n,h}$ (defined in
\eqref{gh})
 as $\tilde g$. Theorem~\ref{Boundedness} says that $\Vert
\widehat{g}_n \Vert_\infty$ is $O_p(1)$, so $\Vert \widehat{f}_n \Vert_\infty$
is also $O_p(1)$, meaning that for any $\epsilon > 0$, there is an $M$ such
that $ \widehat{f}_n$ lies in $\mathcal{LC}_M$ (the class of log-concave
functions bounded by $M$) with probability $1 - \epsilon$.
Thus, we can apply Proposition~\ref{Glivenko1} to conclude that
\begin{equation*}
\left|  \int_\RR \log \lp \widehat{g}_n(t) + b \rp d \lp \GG_n(t) - G^0(t)
  \rp \right|
  \le \sup_{p \in \mathcal{F}_b} \left| \int_\RR p \;  d \lp \GG_n
    -G^0 \rp \right|
  < \epsilon
\end{equation*}
with probability $ 1- \epsilon$ for $n$ large.
Now, let us write
\begin{eqnarray*}
  g^0_h(t) = \piest^0 f^0(\cdot - u^0_1) \star \phi_h
  + (1- \piest^0) f^0(\cdot - u^0_2) \star \phi_h
\end{eqnarray*}
where $\phi_h$ is the density of a $N(0, h^2)$ random variable.
It then
follows from Lemma~\ref{KL2} that as $n \to \infty$,
\begin{equation*}
  \int_\RR \log ( \gestnh) dG^0 \to_p \int_\RR \log(g_h^0) dG^0.
\end{equation*}
Now, since $\int_\RR | \log (g^0) | dG^0 < \infty$,
by Lemma 3.4 of    \cite{seregin_10},
\begin{equation}
  \label{eq:g-consistency-h-b}
  \epsilon(b) + \int \lp \log(g^0(t)+b) - \log(g_h^0(t))\rp dG^0(t)
  \to \int \lp \log(g^0) - \log(g_h^0)\rp dG^0
\end{equation}
as $b \searrow 0$,  by the dominated convergence theorem.
Now, for $|t|$ large enough, $g_h^0(t) \ge g^0(t)$, so $|\log g_h^0(t)| \le |
\log g^0(t)|$, and by \cite{seregin_10}, as was just mentioned, the latter is
$g^0$-integrable, so that
the right side of
\eqref{eq:g-consistency-h-b}
converges to $0$ as $h \searrow 0$.  Thus, for $\epsilon > 0$, we can take
$b$ and $h$ small enough that
\begin{equation*}
  \left|  \epsilon(b)
  + \int_\RR \log(g^0(t)+b) dG^0(t)
  - \int_\RR \log (\gestnh) dG^0 \right| < \epsilon
\end{equation*}
with high probability
for $n$ large enough.  Now, with $h$ fixed,
by  Proposition~\ref{Glivenko2},
\begin{equation*}
  \left| \int_\RR \log(\gestnh) d(\GGn - G^0) \right|
  \le \sup_{p \in \mathcal{H}} \left| \int_\RR p \;  d(\GGn - G^0) \right| < \epsilon
\end{equation*}
with high probability for $n$ large, where $\cal H$ is defined as in
Proposition~\ref{Glivenko2} with $q$ taken as $f_h^0$ and $p_0$ taken as
$g^0$; the proof of Lemma~\ref{KL2} shows that since $(f^0, \phi_h)$ satisfy
Condition $(C)$ then $\log f^0_h$ is $g_0$-integrable, and in fact for any $d
\in \RR$, $\log f^0_h(\cdot -d)$ is $g_0$-integrable (and so is
$f_0$-integrable, so $(f^0, f^0_h)$ satisfy the condition $(C)$ as needed for
Proposition~\ref{Glivenko2}).
 Thus by
\eqref{eq:Hellinger} we are done.  \hfill $\Box$

\bigskip

\par \noindent \textbf{Proof of Corollary~\ref{cor:consistency-f}.} \
 Recall
the definition of $\widecheck{g}_n$ from \eqref{eq:defn:widecheck-g},
as well as
the well-known fact that for any two densities $p_1$ and $p_2$,
\begin{eqnarray}\label{HellingerL1}
  \frac{1}{4} \left(\int_{\RR} \big \vert p_1(t) - p_2(t) \big \vert
    dt\right)^2
  \le H^2(p_1, p_2) \le  \int_{\RR} \big \vert p_1(t) - p_2(t) \big \vert dt.
\end{eqnarray}
Assuming that $\piest^0 < 1/2$ implies by consistency that $\piest < 1/2$
with increasing probability. Hence, the inversion formula (9) in
\cite{bordes2006} yields
\begin{eqnarray}\label{L1L1}
  \int_{\RR} \big \vert \widehat f_n(t) - f^0(t) \big \vert dt & \le &  \frac{1}{1-\piest}  \sum_{k=0}^\infty \left(\frac{\piest}{1-\piest}\right)^k \int_{\RR} \big \vert \widehat g_n(t) - \widecheck{g}_n(t) \big \vert dt  \nonumber \\
  & = & \frac{1}{1-2\piest}   \int_{\RR} \big \vert \widehat g_n(t) - \widecheck{g}_n(t) \big \vert dt.
\end{eqnarray}
From the first inequality in (\ref{HellingerL1}) and Theorem \ref{ConvMixed} above, it follows that
\begin{eqnarray*}
  \int_{\RR} \vert \widehat g_n(t) - g^0(t) \vert dt =o_p(1)
\end{eqnarray*}
which in turn implies that
\begin{eqnarray*}
  \int_{\RR} \vert \widehat g_n(t) - \widecheck{g}_n(t) \vert dt =o_p(1).
\end{eqnarray*}
Indeed,
\begin{eqnarray*}
  \int_{\RR} \vert \widehat g_n(t) - \widecheck{g}_n(t) \vert dt \le \int_{\RR} \vert \widehat g_n(t) - g^0(t) \vert dt + \int_{\RR} \vert g^0(t) - \widecheck{g}_n(t)  \vert dt
\end{eqnarray*}
with
\begin{eqnarray*}
  \int_{\RR} \vert g^0(t) - \widecheck{g}_n(t)  \vert dt &\le &  2 \vert  \piest - \piest^0 \vert  + \int_{\RR} \vert f^0(t - \uestI) - f^0(t - u^0_1) \vert dt \\
  && \ + \int_{\RR} \vert f^0(t - \uestII) - f^0(t - u^0_2) \vert dt \\
  &= & o_p(1)
\end{eqnarray*}
by consistency of $\piest$, $\uestI$ and $\uestII$ and Lemma \ref{Cont}.  It follows now from (\ref{L1L1}) and the inequalities in (\ref{HellingerL1}) that
\begin{eqnarray*}
  H(\widehat f_n, f^0) = o_p(1), \ \textrm{and}  \ \int_{\RR} \vert \widehat f_n(t) - f^0(t) \vert dt = o_p(1).
\end{eqnarray*}
Note that convergence of the MLE $\widehat f_n$ in probability to 0 in the
$L_1$ distance implies its weak convergence to $f^0$ with increasing
probability. Hence, for any arbitrary sequence $(n')$ we can extract a
further subsequence $(n'')$ such that $\widehat f_{n''}$ converges weakly to
$f^0$ almost surely. Hence, the assertion (c) in Proposition 2 of
\cite{cule_10_ejs} holds almost surely for $f_n = \widehat f_{n''}$ and $f =
f^0$, and the remaining claims of our lemma now follow since $(n')$ was
chosen arbitrarily. \hfill $\Box$

\bigskip

Next, we prove the results in Section~\ref{sec:rates}.

\par \noindent \textbf{Proof of Proposition~\ref{Entropy}:} \ It follows from the recent result of \cite{dosswellner13}; see their Theorem 4.1 for $s=0$, that the classes
\begin{eqnarray*}
  \Bigg \{f(\cdot - a), \ f \in \mathcal{SLC}, \  f(0) \in [1/M, M], \ a \in  [a_0 - \delta, a_0+ \delta] \Bigg \}
\end{eqnarray*}
and
\begin{eqnarray*}
  \Bigg \{f(\cdot - b), \ f \in \mathcal{SLC}, \  f(0) \in [1/M, M], \ b \in [b_0 - \delta, b_0+ \delta] \Bigg \}
\end{eqnarray*}
have both the same bracketing entropy $\log(N) \lesssim 1/\sqrt \epsilon$, where $\lesssim$ means that the term on the left side is smaller or equal than the term on the right side up to a positive constant. Let $\{\epsilon, \ldots,  K \epsilon\}$, a $\epsilon$-net for $[0,1]$, where $K = \lfloor 1/\epsilon \rfloor + 1$. For $ i \in \{1, \ldots, N \}$, let $[l_i, u_i]$ and $[l'_i, u'_i]$ an $\epsilon$-bracket for the first and second class respectively. Note that since the brackets are in the Hellinger sense,  we have that $u_i \ge l_i \ge 0 $ and $u'_i \ge l'_i \ge 0$. Now, there exist $i \in \{1, \ldots, N \}$ and $j \in \{1, \ldots, K\}$ such that
\begin{eqnarray*}
  \underline{g} \equiv (j-1) \epsilon l_i + (1-j \epsilon) l'_i \le \lambda f(x- a) + (1- \lambda) f(x- b)  \le \overline{g} \equiv j \epsilon u_i + (1-(j-1) \epsilon) u'_i
\end{eqnarray*}
with
\begin{eqnarray*}
  && H^2\left( \underline{g}, \overline{g}  \right) \\
  && = \frac{1}{2}  \int_{\RR} \left\{ \Big(j \epsilon u_i(t) + (1-(j-1) \epsilon) u'_i(t)\Big)^{1/2} -  \Big((j-1) \epsilon l_i(t) + (1-j \epsilon) l'_i(t)\Big)^{1/2}   \right \}^2 dt  \\
  && \le  2  H^2\Big( (j-1) \epsilon  l_i + (1 - (j-1) \epsilon) u'_i, j \epsilon u_i + (1 - (j-1) \epsilon) u'_i\Big)  \\
  && \ + 2 H^2\Big((1-j \epsilon) l'_i + (j-1) \epsilon l_i, (1-(j-1) \epsilon) u'_i + (j-1) \epsilon l_i\Big)  \\
  && \le 2  H^2\Big( (j-1) \epsilon  l_i, j \epsilon u_i \Big)  +2  H^2\Big((1-j \epsilon) l'_i, (1-(j-1) \epsilon) u'_i\Big), \ \textrm{by Lemma \ref{Hellinger}}  \\
  && \le 4 (j-1) \epsilon \  H^2(l_i, u_i) + 4 \Big(\sqrt{j} - \sqrt{j-1}\Big)^2 \epsilon \\
  && + \  4 (1-j \epsilon) \ H^2(l'_i, u'_i) + 4  \Big( \sqrt{1-(j-1) \epsilon} - \sqrt{1-j \epsilon} \Big)^2,  \ \textrm{applying again Lemma \ref{Hellinger}}.
\end{eqnarray*}
Using the fact that $0 \le j \epsilon \le 1 $ and $1 -(j-1) \epsilon \ge \epsilon$ for all $j \in \{1, \ldots, K \}$, we conclude from the preceding calculations that
\begin{eqnarray*}
  H^2\left( \underline{g}, \overline{g}  \right) \le  4 H^2(l_i, u_i)  +  4 H^2(l'_i, u'_i) +  8 \epsilon  \lesssim \epsilon.
\end{eqnarray*}
The proof is complete by noting that
\begin{eqnarray*}
  \log N_{[ \ ]}(\epsilon,\mathcal G, H)  & \le &  \log K + \log N  \le \log\left(\frac{1}{\epsilon} +1 \right) + \log N\\
  & \le & \frac{1}{\sqrt \epsilon} + \log N \lesssim \frac{1}{\sqrt \epsilon}
\end{eqnarray*}
using the fact that $\log(x +1 ) \le \sqrt x $ for all $x \in [0, \infty)$.   \hfill $\Box$

\bigskip

\bigskip

To prepare for the proof of Theorem \ref{Rates}, we recall that $\piest$, $\uestI$ and $\uestII$ are estimates of $\piest^0$, $u^0_1$ and $u^0_2$ respectively, that are converging at the rate $1/\sqrt n$. Recall also that
$\widecheck{g}_n= \piest f^0(\cdot - \uestI)  + (1 - \piest) f^0(\cdot - \uestII)$, and that $\widehat g_n = \piest \widehat f_n(\cdot - \uestI) + (1-\piest) \widehat f_n(\cdot - \uestII)$, where $\widehat f_n$ is the log-concave MLE. %
As done in \cite{vdvwellner_96} (page 326, Section 3.4), we consider the criterion function
\begin{eqnarray*}
  r_{n,g} = \log \frac{g + \tilde{g}_n}{2 \tilde{g}_n}
\end{eqnarray*}
for $g$ and
\begin{eqnarray*}
  \tilde{g}_n  = \piest \tilde{f}_n(\cdot - \uestI)  + (1- \piest) \tilde{f}_n(\cdot - \uestII)
\end{eqnarray*}
with $\tilde{f}_n \in \mathcal{SCG}_1$ to be constructed. %
Note that
\begin{eqnarray}\label{Remark}
  r_{n, \tilde{g}_n} =0, \textrm{and} \ \mathbb P_n r_{n,\widehat g_n} \ge \frac{1}{2} \PP_n \log \frac{\widehat g_n}{\tilde{g}_n} \ge 0
\end{eqnarray}
where the second claim follows from the definition of the MLE $\widehat g_n$ and concavity of the logarithm. \\

\medskip

Consider now the class of functions
\begin{eqnarray*}
  \mathcal{R}_{n, \eta} & =  & \Big\{r_{n, g} - r_{n, \tilde{g}_n}:g \in \mathcal G, H(g, \tilde{g}_n) < \eta   \Big \} = \Big\{r_{n, g} :g \in \mathcal G, H(g, \tilde{g}_n) < \eta   \Big \}.
\end{eqnarray*}
If $P^0$ denotes again the true probability measure associated with $g^0$,  let
\begin{eqnarray*}
  \mathbb G_n(r_{n, g}) = \sqrt n  \left(\mathbb P_n - P^0\right)r_{n, g} = \sqrt n \left (\frac{1}{n}\sum_{i=1}^n r_{n, g}(X_i) - \int r_{n, g}(x) dP^0(x)\right)
\end{eqnarray*}
and denote
\begin{eqnarray*}
  \Vert \mathbb G_n \Vert_{\mathcal{R}_{n, \eta}} = \sup_{g \in \mathcal R_{n, \eta}} \left \vert \mathbb G_n(r_{n, g}) \right \vert.
\end{eqnarray*}
Finally, define
\begin{eqnarray*}
  \tilde{J}_{[ \ ]}(\delta, \mathcal G, H) = \int_0^\delta \sqrt{1+ \log N_{[ \ ]}(\epsilon,\mathcal G, H)}d\epsilon.
\end{eqnarray*}
Theorem 3.4.4 of \cite{vdvwellner_96} gives sufficient conditions to obtain control on $\Vert \mathbb G_n \Vert_{\mathcal{R}_{n, \eta}}$ in the mean. This control will involve the bracketing entropy bound obtained for the class $\mathcal G$.   One of the crucial conditions to be fulfilled is that the sequence of densities $\tilde{g}_n$ need to be chosen such that $\tilde{g}_n$ approximates the truth $g^0$ and
\begin{eqnarray}\label{Unif}
  \frac{g^0}{\tilde{g}_n} \le M
\end{eqnarray}
for some $M > 0$. Note that the reason we cannot choose $\tilde{g}_n = g^0$ is that in this problem maximization of the log-likelihood involve the random variables $\piest, \uestI$ and $\uestII$, hence it is not at all straightforward to compare the values taken by the criterion at $\widehat g_n$ and $g^0$.  As will be shown in Proposition \ref{tildegn}, we will exhibit an approximating sequence $\tilde g_n$ that will satisfy the condition in (\ref{Unif}) with increasing probability. Based on this proposition, we give now the proof of Theorem \ref{Rates} along the same lines of the proof of Theorem 3.25 in \cite{vdvwellner_96}; see page 290.

\medskip
\medskip

\par \noindent \textbf{Proof of Theorem \ref{Rates}.} \ Fix $\epsilon > 0$.  Let $\tilde{g}_n$ be the approximating sequence of Proposition \ref{tildegn}, and consider the shells
\begin{eqnarray*}
  S_{j,n} = \Big \{g: 2^{j-1} < n^{2/5} H(g, \tilde{g}_n)  \le 2^j \Big \}
\end{eqnarray*}
for integers $j \ge 1$.  Fix an integer $J \ge 1$, and consider the event
\begin{eqnarray*}
  \Big \{ n^{2/5} H(\widehat g_n, \tilde{g}_n) > 2^J \Big \}.
\end{eqnarray*}
Occurrence of this event implies that $\widehat g_n$ belongs to some $S_{j, n}$ with $j \ge J$. But our remark in (\ref{Remark}) implies that
\begin{eqnarray*}
  \sup_{g \in S_{j, n}} \left(\mathbb P_n r_{n, g} - \mathbb P_n r_{n,\tilde{g}_n} \right)   = \sup_{g \in S_{j, n}} \mathbb P_n r_{n, g} \ge 0.
\end{eqnarray*}
Thus,  for any $\delta > 0$ and $M > 0$ we can write
\begin{eqnarray*}
  P\left( n^{2/5} H(\widehat g_n, \tilde{g}_n) > 2^J \right)  & = & P\left( n^{2/5} H(\widehat g_n, \tilde{g}_n) > 2^J, g^0 \le M \tilde{g}_n  \right)  + P\left(g^0 > M \tilde{g}_n \right)  \\
  & = & \sum_{j \ge J, 2^j \le n^{2/5} \eta} P\left( \sup_{g \in S_{j, n}} \left(\mathbb P_n r_{n, g} - \mathbb P_n r_{n,\tilde{g}_n} \right) \ge 0, g^0 \le M \tilde{g}_n \right) \\
  && + \  P\left( H(\widehat g_n, \tilde{g}_n)  \ge \delta/2 \right)  + P\left(g^0 > M \tilde{g}_n \right).
\end{eqnarray*}
By Proposition \ref{tildegn}, for any $\eta > 0$ and $M$ large enough the second and third terms are bounded by $\epsilon/2$.  Now, using Theorem 3.4.4 of \cite{vdvwellner_96}, we have
\begin{eqnarray}\label{Diff}
  P^0(r_{n, g} - r_{n, \tilde{g}_n}) \lesssim -H^2(g, \tilde{g}_n)
\end{eqnarray}
for all densities $g$ such that $H(g, \tilde{g}_n) \ge 32 M H(\tilde{g}_n, g^0)$, and
\begin{eqnarray}\label{Control}
  E_{P^0} \Vert \mathbb G_n \Vert_{\mathcal{R}_{n, \eta}}  \lesssim \tilde{J}_{[\ ]}(\eta, \mathcal{G}, H) \left(1+ \frac{\tilde{J}_{[ \ ]}(\eta, \mathcal{G}, H)}{\eta^2 \sqrt n}  \right).
\end{eqnarray}
From Proposition \ref{Entropy}, we know that there exists a constant $K > 0$ (not depending on $n$) such that
\begin{eqnarray*}
  \log N_{[ \ ]}(\epsilon,\mathcal G, H)   \le \frac{K}{\sqrt{\epsilon}}.
\end{eqnarray*}
For $\eta > 0$ small enough, $1 < K/\sqrt{\epsilon}$ and hence
\begin{eqnarray*}
  \tilde{J}_{[\ ]}(\eta, \mathcal{G}, H)  &  \le &  \int_0^\eta \frac{\sqrt{2K}}{\epsilon^{1/4}} d\epsilon = \frac{4 \sqrt{2K}}{3} \eta^{3/4}.
\end{eqnarray*}
Now, define
\begin{eqnarray*}
  \phi_n(\eta):= \frac{4 \sqrt{2K}}{3} \eta^{3/4}  \left(1 + \frac{4 \sqrt{2K}}{3 \eta^{5/4} \sqrt n} \right).
\end{eqnarray*}
We have that
\begin{eqnarray*}
  \frac{\phi_n(\eta)}{\eta} \propto  \frac{1}{\eta^{1/4}} \left(1 + \frac{4 \sqrt{2K}}{3 \eta^{5/4} \sqrt n} \right)
\end{eqnarray*}
and hence the function $\eta \mapsto \phi_n(\eta)/\eta$ is decreasing on $(0, \infty)$. Also, if we put $r_n  = n^{2/5}$ then,
\begin{eqnarray*}
  r_n^2 \phi_n\left( \frac{1}{r_n} \right)  &=  & n^{4/5} \frac{4 \sqrt{2K}}{3} n^{-3/10} \left(1 + \frac{4 \sqrt{2K}}{3} \right) \\
  & =& \frac{4 \sqrt{2K}}{3} \left(1 + \frac{4 \sqrt{2K}}{3} \right) \ \sqrt n.
\end{eqnarray*}
Now note that on the event $\Big \{g^0 \le M \tilde{g}_n \Big  \} $ it follows from (\ref{Diff}) that for all $g \in S_{n,j}$
\begin{eqnarray*}
  P^0(r_{n, g} - r_{n, \tilde{g}_n}) \le \frac{-2^{2j -2}}{r^2_n}
\end{eqnarray*}
and hence %
\begin{eqnarray}\label{ResRates}
  &&\sum_{j \ge J, 2^j \le n^{2/5} \eta} P\left( \sup_{g \in S_{j, n}} \left(\mathbb P_n r_{n, g} - \mathbb P_n r_{n,\tilde{g}_n} \right) \ge 0, g^0 \le M \tilde{g}_n \right) \nonumber \\
  && \le \sum_{j \ge J, 2^j \le n^{2/5} \eta} P\left( \Vert \mathbb G_n \Vert_{S_{n,j}}   \ge \frac{-2^{2j -2} \sqrt n}{r^2_n} \right)  \nonumber \\
  && \le \sum_{j \ge J} \frac{\phi_n(2^j/r_n) r^2_n}{2^{2j -2} \sqrt n}, \ \textrm{by Markov's inequality} \nonumber \\
  &&  \le \frac{1}{4} \sum_{j \ge J} \frac{\phi_n(1/r_n) r^2_n}{2^j \sqrt n},  \ \textrm{using the fact that} \frac{\phi_n(ct)}{ct} \le \frac{\phi_n(t)}{t}  \textrm{for $c > 1$} \nonumber \\
  && \lesssim \sum_{j \ge J} 2^{-j}  \to 0 \textrm{ as $J \to \infty$}
\end{eqnarray}
which implies that  $H(\widehat g_n, \tilde{g}_n) = O_p(n^{-2/5})$. Hence, $L_1(\widehat g_n, \tilde{g}_n) = O_p(n^{-2/5})$ which in turn implies that $L_1(\widehat{g}_n, g^0) = O_p(n^{-2/5})$ using the result of  Proposition \ref{tildegn} and the triangular inequality.

To conclude a similar rate result for $\widehat f_n$, we show now that $L_1(\widehat f_n, c_n \tilde{f}_n) = O_p(n^{-2/5})$ where $\tilde{f}_n$ and the normalizing constant $c_n$ are defined in Proposition \ref{tildegn}. Using again the inversion formula (9) in \cite{bordes2006}, we can write
\begin{eqnarray*}
  \int_{\RR} \big \vert \widehat{f}_n(t) - c_n \tilde{f}_n(t) \big \vert dt \le \frac{1}{1-2\piest}\int_{\RR} \big \vert \widehat{g}_n(t) - \tilde{g}_n(t) \big \vert dt  = O_p(n^{-2/5}).
\end{eqnarray*}
The proof is complete using the result of Proposition \ref{tildegn} and the
triangle inequality.  \hfill $\Box$

\medskip

\begin{prop}\label{tildegn}
  There exist $\tilde{f}_n \in \mathcal{SCG}$ and $c_n > 0$ such that $c_n \tilde{f}_n \in \mathcal{SCG}_1$ and
  \begin{eqnarray*}
    \sup_{t \in \RR} \frac{g^0(t)}{\tilde{g}_n(t)} = O_p(1)
  \end{eqnarray*}
  where
  \begin{eqnarray*}
    \tilde{g}_n = \piest c_n\tilde{f}_n(\cdot  - \uestI) + (1- \piest) c_n\tilde{f}_n(\cdot - \uestII).
  \end{eqnarray*}
  Furthermore, if $\sqrt n (\uestI - \uestI^0) = O_p(1), \sqrt n (\uestII - u^0_2) =
  O_p(1)$ and $\sqrt{n} (\piest - \piest^0) = O_p(1/\sqrt{n})$, then
  \begin{eqnarray*}
    L_1(\tilde{f}_n, f^0) = O_p(n^{-1/2}), \ \ \textrm{and} \ \ L_1(\tilde{g}_n, g^0) = O_p(n^{-1/2}).
  \end{eqnarray*}

\end{prop}

\medskip

\par \noindent \textbf{Proof of Proposition \ref{tildegn}.} \ In the following we denote
\begin{eqnarray*}
  \delta_1 = \uestI - u^0_1 \ \ \textrm{and} \ \ \delta_2 = \uestII - u^0_2.
\end{eqnarray*}
Suppose that $0 \le \delta_2 \le \delta_1$. Define
\begin{eqnarray*}
  \tilde{f}_n(t)  = \left \{
    \begin{array}{lll}
      f^0(\vert t \vert - \delta_1), \ \ \textrm{if $t \le - \delta_1$} \\
      f^0(0), \hspace{1.2cm} \textrm{if $ -\delta_1 \le t \le \delta_1$} \\
      f^0(t  - \delta_1), \ \ \ \ \textrm{if $t \ge  \delta_1$}.
    \end{array}
  \right.
\end{eqnarray*}
It is not difficult to see that $\tilde{f}_n \in \mathcal{SLC}$. Now, define the ratios
\begin{eqnarray*}
  \tilde{R}_{n,1}(t): = \frac{f^0(t)}{\tilde{f}_n(t-\delta_1)}, \ \ \textrm{and} \ \ \tilde{R}_{n,2}(t): = \frac{f^0(t)}{\tilde{f}_n(t-\delta_2)}.
\end{eqnarray*}
We have that
\begin{eqnarray*}
  \tilde{R}_{n,1}(t) = \left \{
    \begin{array}{lll}
      \frac{f^0(t)}{f^0(\delta_1 - t - \delta_1)} = 1, \ \ \textrm{if $t \le 0$} \\
      \frac{f^0(t)}{f^0(0)}  \le 1, \hspace{1.2cm} \textrm{if $ 0 \le  t \le 2 \delta_1  $ } \\
      \frac{f^0(t)}{f^0(t - 2 \delta_1 ) } \le 1, \ \ \ \ \ \ \textrm{if $t \ge 2 \delta_1$ }.
    \end{array}
  \right.
\end{eqnarray*}
The third inequality follows from the fact that $f^0$ is decreasing on $[0, \infty)$.  Also,
\begin{eqnarray*}
  \tilde{R}_{n,2}(t) = \left \{
    \begin{array}{lll}
      \frac{f^0(t)}{f^0(\delta_2 - \delta_1 - t)} \le 1 , \ \ \textrm{if $t \le \delta_2 - \delta_1$} \\
      \frac{f^0(t)}{f^0(0)}  \le 1, \hspace{1.3cm} \textrm{if $ \delta_2 - \delta_1 \le  t \le \delta_1  + \delta_2 $ } \\
      \frac{f^0(t)}{f^0(t -  \delta_1 - \delta_2)} \le 1, \ \ \  \textrm{if $t \ge  \delta_1 + \delta_2$}.
    \end{array}
  \right.
\end{eqnarray*}
The first inequality follows from symmetry of $f^0$ which allows to write that $f^0(t)/f^0(\delta_2 - \delta_1 - t)  = f^0(x)/f^0(x  - (\delta_1 - \delta_2))$ with $x = -t \ge \delta_1 - \delta_2 \ge 0$.  Since $f^0$ is decreasing on $[0, \infty)$, it follows that $f^0(x) \le f^0(x  - (\delta_1 - \delta_2))$. The third inequality is again a consequence of the latter property of $f^0$.

Now let $c_n = \left(\int_{\RR} \tilde{f}_n(t) dt\right)^{-1}$. We have that
\begin{eqnarray*}
  c^{-1}_n -1 \  & =  &  \int_{-\infty}^{-\delta_1}  f^0(\vert t \vert - \delta_1) dt + \int_{\delta_1} ^\infty f^0(t - \delta_1) dt + 2 \delta_1 f^0(0) - 1  \\
  & = & 2 \delta_1 f^0(0).
\end{eqnarray*}
This in turn implies that
\begin{eqnarray*}
  0 \le 1-c_n = \frac{2 \delta_1 f^0(0)}{1+ 2 \delta_1 f^0(0)}.
\end{eqnarray*}
Now, define $\tilde{g}_n:= c_n \left(\piest \tilde{f}_n(\cdot - \uestI) + (1-\piest) \tilde{f}_n(\cdot - \uestII)\right)$. Then,
\begin{eqnarray*}
  \frac{g^0(t)}{\tilde{g}_n(t)} & = & c^{-1}_n \frac{\piest^0 f^0(t - u^0_1) +(1-\piest^0) f^0(t - u^0_2)}{\piest \tilde{f}_n(t - \uestI) + (1-\piest) \tilde{f}_n(t- \uestII)} \\
  & \le & c^{-1}_n \frac{\piest^0}{\piest} \frac{f^0(t - u^0_1)}{\tilde{f}_n(t - u^0_1 - \delta_1)} + c^{-1}_n \frac{1-\piest^0}{1-\piest} \frac{f^0(t -u^0_2)}{\tilde{f}_n(t - u^0_2 - \delta_2)} \\
  & = & c^{-1}_n \frac{\piest^0}{\piest} \tilde{R}_{n,1}(t-u^0_1) + c^{-1}_n \frac{1-\piest^0}{1-\piest}  \tilde{R}_{n,2}(t- u^0_2) \\
  & \le& c^{-1}_n \left(\frac{\piest^0}{\piest} + \frac{1-\piest^0}{1-\piest}\right).
\end{eqnarray*}
By consistency of $\piest$ and $\uestII$, we have $ \piest^0 < \piest \le (3/2) \piest^0$ and $c^{-1}_n \le 3/2$ with increasing probability. Hence, we can bound the right hand side of the preceding display by $9/2$ with increasing probability. As the other cases can be handled similarly, the details are skipped but we give below the corresponding expression of $\tilde{f}_n$:

\begin{itemize}
\item If $0 \le \delta_1 < \delta_2$, then  we only need to switch the roles of $\delta_1$ and $\delta_2$ and hence take
  \begin{eqnarray*}
    \tilde{f}_n(t)  = \left \{
      \begin{array}{lll}
        f^0(\vert t \vert - \delta_2), \ \ \textrm{if $t \le - \delta_2$} \\
        f^0(0), \hspace{1.2cm} \textrm{if $ -\delta_2 \le t \le \delta_2$} \\
        f^0(t  - \delta_2), \ \ \ \ \textrm{if $t \ge  \delta_2$}.
      \end{array}
    \right.
  \end{eqnarray*}

\item If $\delta_2 \le  0 \le \delta_1$, then we can take
  \begin{eqnarray*}
    \tilde{f}_n(t)  = \left \{
      \begin{array}{lll}
        f^0(\vert t \vert - (\delta_1 - \delta_2)), \ \ \textrm{if $t \le - (\delta_1 - \delta_2)$} \\
        f^0(0), \hspace{1.2cm} \textrm{if $ -(\delta_1 - \delta_2) \le t \le \delta_1 - \delta_2$} \\
        f^0(t  -(\delta_1 - \delta_2)), \ \ \ \ \textrm{if $t \ge  \delta_1 - \delta_2$}.
      \end{array}
    \right.
  \end{eqnarray*}

\item If $\delta_1 \le 0 \le\delta_2 \le 0$,  then we only need to switch $\delta_1$ and $\delta_2$.

\item If $\delta_1 \le \delta_2 \le 0$, then we can take

  \begin{eqnarray*}
    \tilde{f}_n(t)  = \left \{
      \begin{array}{lll}
        f^0(\vert t \vert  + \delta_1 ), \ \ \textrm{if $t \le \delta_1 $} \\
        f^0(0), \hspace{1.2cm} \textrm{if $  \delta_1 \le t \le -\delta_1 $} \\
        f^0(t  +\delta_1 ), \ \ \ \ \textrm{if $t \ge  -\delta_1 $}.
      \end{array}
    \right.
  \end{eqnarray*}

\item If $\delta_2 \le \delta_1 \le 0$, we again switch the roles of $\delta_1$ and $\delta_2$.
\end{itemize}

\medskip

In all the cases above, one can verify that the ratios $\tilde{R}_{n,1}$ and $\tilde{R}_{n,2}$ as defined above stay below $1$. We would like to stress the fact that the way $\tilde{f}_n$ is constructed is not unique: one only need to exhibit examples which would give control of the ratio $g^0/\tilde{g}_n$. To show now the second assertion,   we will again consider only the first case where $0 \le \delta_2 \le \delta_1$ since the remaining configurations can be handled similarly.   We have that
\begin{eqnarray*}
  L_1(\tilde{f}_n, f^0) & = & \int_{\RR} \vert \tilde{f}_n(t) - f^0(t)  \vert dt \\
  & = & \int_{-\infty}^{-\delta_1} \vert f^0( - t - \delta_1) - f^0(t) \vert dt + \int_{-\delta_1}^{\delta_1} \vert f^0(0) - f(t) \vert dt \\\
  && \ + \int_{\delta_1}^{\infty} \vert f^0( t - \delta_1) - f^0(t) \vert dt \\
  & = & 2 \int_{\delta_1}^\infty \vert f^0(t - \delta_1) - f^0(t) \vert dt + 2 \delta_1 f^0(0)  \\
   &  \le & 2 (C+1) \delta_1
\end{eqnarray*}
where $ C$ is the constant given in Lemma \ref{Cont}.  By the assumption on the rate of convergence of $\delta_1$, this implies that $L_1(\tilde{f}_n, f^0) = O_p(n^{-1/2})$. Also, we have
\begin{eqnarray*}
  && L_1(\tilde{g}_n, g^0) \\
  & &= \int_{\RR} \vert \tilde{g}_n(t) - g^0(t) \vert dt  \\
  && \int_{\RR} \left \vert c_n \left(\piest \tilde{f}_n(t-\uestI) + (1-\piest) \tilde{f}_n(t-\uestII)   \right)  - \left( \piest^0 f^0(t - \uestI^0) + (1-\piest^0) f^0(t - u^0_2 )\right)  \right\vert dt \\
&& \le 1-c_n + \int_{\RR} \left \vert \piest \tilde{f}_n(t-\uestI)  - \piest^0 f^0(t - \uestI^0)  \right \vert dt \\
&& + \int_{\RR} \left \vert (1-\piest) \tilde{f}_n(t-\uestII)   - (1-\piest^0) f^0(t - u^0_2 ) \right\vert dt \\
&& = \le 1-c_n + 2  \vert \piest - \piest^0 \vert  + \int_{\RR} \left \vert \tilde{f}_n(t-\uestI)  - f^0(t-u^0_1) \right \vert dt + \int_{\RR} \left \vert \tilde{f}_n(t-\uestII)  - f^0(t-u^0_2) \right \vert dt \\
&& =  \le 1-c_n + 2  \vert \piest - \piest^0 \vert  + \int_{\RR} \left \vert \tilde{f}_n(t-\delta_1) - f^0(t) \right \vert dt + \int_{\RR} \left \vert \tilde{f}_n(t-\delta_2)  - f^0(t) \right \vert dt.
\end{eqnarray*}
Using the definition of $\tilde{f}_n$ for the first case, $0 \le \delta_2 < \delta_1$, we can write
\begin{eqnarray*}
\int_{\RR} \left \vert \tilde{f}_n(t-\delta_1) - f^0(t) \right \vert dt   &= & \int_{-\infty}^0 \left \vert f^0(\delta_1 -t -\delta_1) - f^0(t) \right \vert + \int_0^{2\delta_1} \left \vert f^0(0) - f^0(t) \right \vert dt \\
&& \ + \int_{2\delta_1}^\infty \left \vert f^0(t-2\delta_1) - f^0(t) \right \vert dt \\
& \le  &2 \delta_1 f^0(0)  + 6 \delta_1  f^0(0), \ \textrm{using Lemma A.1}\\
& = & 8 \delta_1 f^0(0) \\
  && = O_p(n^{-1/2})
\end{eqnarray*}
using again Lemma \ref{Cont},  and the assumption on the rate of convergence of $\piest$ and $\uestI$ and $\uestII$.

\hfill $\Box$

\bigskip

\subsection{Auxiliary Results}

\begin{lemma}\label{Cont}
  Let $f \in \mathcal{SLC}_1$. Then, there exists a constant $C >0$ depending only on $f$ such that
  \begin{eqnarray*}
    \int_{\RR} \Big \vert f(t + \delta) - f(t) \Big \vert dt \le C \delta
  \end{eqnarray*}
  for all $\vert \delta \vert \le 1$.
\end{lemma}

\medskip

\par \noindent \textbf{Proof.} \ For $\delta \ge 0$, we have $f(t+\delta) \le f(t)$ on $t \in [0, \infty)$, $f(t - \delta) \ge f(t)$ on $t \in [\delta/2, \infty)$ and $f(t-\delta) \le f(t)$ on $[0, \delta/2]$. Using  the symmetry of $f$, we can write
\begin{eqnarray*}
  \int_{\RR} \Big \vert f(t+ \delta) - f(t) \Big \vert dt &= & \int_0^\infty \Big \vert f(t + \delta) - f(t) \Big \vert  dt+ \int_{0}^\infty \Big \vert f(t -\delta) - f(t) \Big \vert dt \\
  & = & \int_0^\infty \Big(f(t) - f(t+\delta) \Big) dt + \int_0^{\delta/2}\Big( f(t) -f(t - \delta) \Big)dt \\
  && \ + \int_{\delta/2}^\infty \Big (f(t - \delta) - f(t) \Big)dt  \\
  & \le & 1/2 - \int_{\delta}^\infty f(t) dt  +  \int_0^{\delta/2} f(t) dt + \int_{-\delta/2}^\infty f(t) dt - \int_{\delta/2}^\infty  f(t) dt \\
  & = & \int_0^\delta f(t) dt  + \int_0^{\delta/2} f(t) dt   + 2 \int_0^{\delta/2} f(t) dt \\
  & \le & 3 \delta \sup_{t \in [0,1]} f(t) = 3 \delta f(0)
\end{eqnarray*}
provided that $\delta \le 1$ . The same reasoning can be applied for negative
values of $\delta$. \hfill $\Box$

\bigskip
\bigskip

\begin{lemma}
  \label{Hellinger}
  For any positive functions $p, q, h$ we have that
  \begin{eqnarray*}
    H(p + h, q + h) \le H(p, q)
  \end{eqnarray*}
\end{lemma}

\medskip

\par \noindent \textbf{Proof.} \ By definition, we have that
\begin{eqnarray*}
  2 H^2(p+h, q+h) & =  & \int_{\RR}  \left( \sqrt{p(t) + h(t)} - \sqrt{q(t) + h(t)} \right)^2 dt \\
  & = & \int_{\RR} \left(\frac{p(t)  - q(t)}{\sqrt{p(t) + q(t)} + \sqrt{q(t) + h(t)}} \right)^2 dt \\
  & = & \int_{\RR} \left(\sqrt{p(t)} - \sqrt{q(t)}\right)^2   \left(\frac{\sqrt{p(t)} + \sqrt{q(t)}}{\sqrt{p(t) + h(t)} + \sqrt{q(t) + h(t)}}\right)^2 dt \\
  & \le &  \int_{\RR} \left(\sqrt{p(t)} - \sqrt{q(t)}\right)^2 dt = 2 H^2(p, q)
\end{eqnarray*}
where the last inequality follows since $h \ge 0$, hence the result.  \hfill $\Box$

\medskip
\medskip

\par \noindent \textbf{The function $J$ and its partial derivatives}: As in
\cite{duembgen_09}, we consider the two-dimensional function $J$ defined by
\begin{eqnarray*}
  J(r, s) = \int_0^1 \exp((1-t) r + ts ) dt
\end{eqnarray*}
for $r, s \in \RR$.  Using the same notation of these authors, define
\begin{eqnarray*}
  J_{a,b}(r, s) = \frac{\partial^{a+b} J(r, s)}{\partial^a s \partial^b s}.
\end{eqnarray*}
Direct calculations yield
\begin{eqnarray*}
  J_{a, b}(r, s) = \exp(r) J_{a, b}(0, s-r)
\end{eqnarray*}
with
\begin{eqnarray*}
  J_{0,0}(0, y) = J(0, y) =
  \begin{cases}
    1 , \ \textrm{if $y=0$} \\
    \frac{\exp(y)-1}{y}, \ \textrm{otherwise},
  \end{cases}
\end{eqnarray*}
\begin{eqnarray*}
  J_{0,1}(0, y) =
  \begin{cases}
    \frac{1}{2} , \ \textrm{if $y=0$} \\
    \frac{y\exp(y)-(\exp(y) -1)}{y^2}, \ \textrm{otherwise},
  \end{cases}
\end{eqnarray*}
and
\begin{eqnarray*}
  J_{0,2}(0, y)
  \begin{cases}
    \frac{1}{3}, \ \textrm{if $y=0$} \\
    \frac{y^2\exp(y)-2y \exp(y) + 2(\exp(y) -1)}{y^3}, \ \textrm{otherwise}.
  \end{cases}
\end{eqnarray*}

\bibliographystyle{ims}

\begin{all-in-one-file}
  \bibliography{MixSym}
\end{all-in-one-file}

\end{document}